\documentclass[twoside,11pt]{article}

\usepackage{jmlr2e}
\usepackage{amsmath}
\usepackage{amssymb}

\usepackage{float}

\usepackage{amsopn}
\usepackage{graphicx}
\usepackage{latexsym}
\usepackage{stmaryrd}
\usepackage{algorithm}
\usepackage{algorithmic} 
\usepackage{enumerate}
\usepackage{natbib}
\usepackage{indentfirst}

\ShortHeadings{Learning the Sparse and  Low Rank PARAFAC Decomposition via the Elastic Net}{}

\firstpageno{1}

\begin{document}

\title{Learning the Sparse and  Low Rank PARAFAC Decomposition via the Elastic Net}
\author{\name }
\author{\name Songting Shi \email songtingshists@pku.edu.cn \\
       \addr Department of Scientific and Engineering Computing\\
       School of Mathematical Sciences\\
      	Peking University\\
	Beijing 300071, P. R. China\\
       \AND
       \name Xiang Li \email xli60@mgh.havard.edu \\
       \addr Clinical Data Science Center\\
       Massachusetts General Hospital\\
       Boston, MA 02114, USA\\
       \AND
      \name Arkadiusz Sitek \email sarkadiu@gmail.com \\
       \addr Gordon Center for Medical Imaging\\
       Massachusetts General Hospital\\
       Boston, MA 02114, USA\\
       \AND
       \name Quanzheng Li \email li.quanzheng@mgh.harvard.edu \\
       \addr Gordon Center for Medical Imaging\\
       Massachusetts General Hospital\\
       Boston, MA 02114, USA}

%\author{\name Songting Shi \email songtingshists@pku.edu.cn \\
%       \addr Department of Scientific and Engineering Computing\\
%       School of Mathematical Sciences\\
%      	Peking University\\
%	Beijing 300071, P. R. China\\
%       \AND
%       \name Xiang Li \email xli60@mgh.havard.edu \\
%       \addr Division of Computer Science and Department of Statistics\\
%       University of California\\
%       Berkeley, CA 94720-1776, USA\\
%       \AND
%      \name Arkadiusz Sitek \email sarkadiu@gmail.com \\
%       \addr Division of Computer Science and Department of Statistics\\
%       University of California\\
%       Berkeley, CA 94720-1776, USA
%       \AND
%       \name Quanzheng Li \email li.quanzheng@mgh.harvard.edu \\
%       \addr Division of Computer Science and Department of Statistics\\
%       University of California\\
%       Berkeley, CA 94720-1776, USA}

%\editor{Kevin Murphy and Bernhard Sch{\"o}lkopf}
%\editor{}
\maketitle

\begin{abstract}%   <- trailing '%' for backward compatibility of .sty file
In this article, we derive  a Bayesian model to learning the sparse and low rank PARAFAC decomposition for the observed tensor with missing values via the elastic net, with property to find the true rank and  sparse factor matrix which is robust to the noise. We formulate  efficient  block coordinate descent  algorithm and  admax stochastic block coordinate descent  algorithm to solve it, which can be used to solve the large scale problem. To choose the appropriate rank and sparsity in PARAFAC decomposition, we will give a solution path by gradually increasing the regularization to increase the sparsity and decrease the rank. When we find the sparse structure of the factor matrix, we can fixed the sparse structure, using a small to regularization to decreasing the recovery error, and one can choose the proper decomposition from the solution path with sufficient sparse factor matrix with low recovery error.  We test the power of  our algorithm on the simulation data and real data, which show it is powerful.
\end{abstract}

\begin{keywords}
  PARAFAC decomposition, tensor imputation, elastic net
\end{keywords}

\section{Introduction}
Tensor decomposition has been used as a powerful tool on many fields (\cite{td_review}) for multiway data analysis. Among the representation models for tensor decompositions, PARAFAC decomposition is one of the simplest form. It aims to extract the data-dependent basis functions (i.e. factors), so that the observed tensor can be expressed by the multi-linear combination of these factors, which can reveal meaningful underlying structure of the data. The ideal situation for the PARAFAC decomposition is that we know (1) the true rank of the tensor, and (2) the sparse structure (i.e. locations of zeros) of the factors. Given a wrong estimation of the tensor's rank or sparse structure of the factors, even the recovery error of the decomposition results are very low, there still exits a high generalization error as the learning procedure tries to fit the noise information. But in most of the real-case scenarios, we won't be able to obtain the prior knowledge of the true rank of the tensor and the true sparse structure of the factors. Various literatures have proposed solutions for either of the above problems. For example, Liu et al. proposed the sparse non-negative tensor factorization using columnwise coordinate descent to get the sparse factors
%\yrcite{sntf_ccd}.
(\cite{sntf_ccd}).
 To get the proper rank estimation, a novel rank regularization PARAFAC decomposition was proposed in (\cite{rankreg}). However, to the author's best knowledge, there is no work that can handle both of the two problems simultaneously  in an integrated framework. And in many works, they  try to emphasize on the lower recovery error, few authors focus on finding the true underline factors of the observed data with noise. As we can see from our simulation examples, there are algorithms( e.g. CP-ALS), although they find a lower relative recovery error, but they deviate from the true underline factors. The reason is that these algorithm learn  to fit the noise information to reduce the recovery error. Our algorithm is effective  to de-noise from both our simulation data and real dataset. Moreover, there seems few researches to use stochastic scheme to calculate  the large scale PARAFAC decomposition.

In response to the above problems, we develop a PARAFAC decomposition model capturing the true rank and producing the sparse factors by minimizing the least square lost with elastic net regularization.

The elastic net regularization is a convex combination of the $l_2$ regularization and $l_1$ regularization, where the $l_2$ regularization
helps to capture the true rank, 
while the $l_1$ regularization trends to yield the sparse factors. 
For solving the minimization problem, we use the block coordinate descent algorithm and stochastic block coordinate descent algorithm, where the block coordinate descent algorithm 
  have been shown to be an efficient scheme (\cite{bcd}) especially on large scale problems.
To identify the true rank and sparse factors, we generate a solution path  by gradually increasing the magnitude of the regularization. And the initial dense decomposition is gradually transformed to a low rank and sparse decomposition. 
But when we find the  rank and the appropriate sparse factors, the recovery error may be  large because large regularization will shrink  the solution closer to origin.
 One can fix the sparse structure of the factors and perform a constrained optimization with small regularization to decrease the recovery error. 
 At last, we can choose the appropriate low rank and sparse PARAFAC decomposition from  the refined solution of the solution path which features a low recovery error. 
  We test our algorithm firstly on the synthetic data. The results shows 
 that the algorithm captures the true rank with high probability and the sparse structure of the factors found by the algorithm is sufficiently close to the sparse structure of the true factors,  and  the sparse constraint algorithm with the sparse structure of the factors  trends to find the true underline factors.
And it finds out the meaningful structure of features when applied on the coil-20 dataset and  the COPD data.

%The rest of this paper is organized as follows: We  derive our model from the Bayesian framework in Section \ref{sec:bays_infer}; The coordinate descent algorithm to solve the model
%% The main algorithm which use the coordinate descent to achieve the solution of the model 
%  is proposed in Section \ref{sec:main_algorithm}; the theory analysis is given in Section \ref{sec:theory}; the solution path method in practical is given in Section \ref{sec:solution_path}; report the numerical test in Section \ref{sec:numerical_test};  and conclude this paper in Section \ref{sec:conclusion}; 

The rest of this paper is organized as follows: 
%We derive our model from the Bayesian framework in Section \ref{sec:bays_infer}. 
Our model based on a Bayesian framework is derived in Section  \ref{sec:bays_infer}. The coordinate descent algorithm to solve the optimization is proposed in Section \ref{sec:main_algorithm}. Theoretical analysis is given in Section \ref{sec:theory}. The solution path method is given in Section \ref{sec:solution_path}.  The   stochastic block coordinate descent algorithm to solve the optimization is proposed in Section \ref{sec:main_algorithm_random}. Results from applying the algorithm on the synthetic data and real data is summarized in Section \ref{sec:numerical_test}, and conclude this paper in Section \ref{sec:conclusion}.

%\section{Problem Formulation}
\section{Bayesian PARAFAC model}
\label{sec:bays_infer}
The notations used through out the article are as follows: We use bold lowercase and capital letters for vectors ${\bf a}$, and matrices ${\bf A}$, respectively. Tensors are underlined, $\text{e.g.}$, $\underline{\bf X}$. Both the matrix and tensor Frobenius norms are represented by $||\cdot||_ F$ and
% use $|\cdot |$ denote the $l_1$ norm f. 
use $\text{vec}(\cdot)$ denote the vectorization of a matrix.
  Symbols $\varotimes$, $\varodot$, $\varoast$, and $\circ$, denote the Kroneker, Kathri-Rao, Hadamard (entry-wise), and the outer product, respectively. $\varoslash$ denote the entry-wise divide operation, and  $[N] := \{1, 2, \ldots, N\}$. $\text{Tr}(\cdot)$ stands for the trace operation.

The observation $\underline{\bf  Z } \in \mathbb R^{I_1\times I_2 \times \cdots \times I_N}$  is composed of the true solution $\underline{\bf  X } \in \mathbb R^{I_1\times I_2 \times \cdots \times I_N}$ and the noise tensor $\underline {\bf E} \in \mathbb  R^{I_1\times I_2 \times \cdots \times I_N}$.
\begin{equation}\label{Z_data}
 \underline{\bf  Z }_{i_1, i_2, \ldots, i_N} = \underline{\bf  X }_{i_1, i_2, \ldots, i_N} +\underline{\bf  E }_{i_1, i_2, \ldots, i_N}
 \end{equation}
 where $ \underline{\bf  E }_{i_1, i_2, \ldots, i_N} \sim \mathcal N(0, \sigma^2), \; \text{i.i.d.}$.
  We only observed a subset entries of  $\underline {\bf Z}$, which is given by a binary tensor $\underline {\boldsymbol {\Delta}} \in \{0, 1\} ^{I_1\times I_2 \times \cdots \times I_N}$.
%  , where $1$ stands for the location of known values, $0$ i.

%Now we begin to introduce our model. 
%The  observation $\underline{\bf  Z } \in \mathbb R^{I_1\times I_2 \times \cdots \times I_N}$ is composed of the true solution $\underline{\bf  X } \in \mathbb R^{I_1\times I_2 \times \cdots \times I_N}$ and the noise tensor $\underline {\bf E} \in \mathbb  R^{I_1\times I_2 \times \cdots \times I_N}$.
%\begin{equation}\label{Z_data}
% \underline{\bf  Z }_{i_1, i_2, \ldots, i_N} = \underline{\bf  X }_{i_1, i_2, \ldots, i_N} +\underline{\bf  E }_{i_1, i_2, \ldots, i_N}
% \end{equation}
% where $ \underline{\bf  E }_{i_1, i_2, \ldots, i_N} \sim \mathcal N(0, \sigma^2), \; \text{i.i.d.}$.
%  And we only observed a subset entries of  $\underline {\bf Z}$, which is given by a binary tensor $\underline {\boldsymbol {\Delta}} \in \{0, 1\} ^{I_1\times I_2 \times \cdots \times I_N}$.
%  , where $1$ stands for the location of known values, $0$ i.

Assume that $\underline{\bf  X }$ comes from the PARAFAC decomposition $\underline{\bf  X } = {\bf A}^{(1)} \circ {\bf A}^{(2)}  \circ \cdots \circ {\bf A}^{(N)}  $, which means that
\begin{equation}
\underline{\bf  X }_{i_1, i_2, \ldots, i_N}   = \sum_{r=1}^R {\bf A}^{(1)} (i_1,r)  {\bf A}^{(2)} (i_2,r)    \cdots  {\bf A}^{(N)} (i_N,r)  
\end{equation}
where the factor matrix ${\bf A }^{(n)}  \in \mathbb{R}^{I_n \times R}$,  for $n = 1, 2, \ldots, N$.
 Let ${\bf a}^{(n)}_r, \; r = 1, \cdots R $ be the columns of ${\bf A }^{(n)}$.
%  Let ${\bf a}_r^{\circ} := {\bf a}_r^{(1)} \circ \cdots \circ {\bf a}_r^{(N)}$
 Then we may express the decomposition  as 
 \begin{equation}\label{cp_decompositon}
  \begin{aligned}
\underline{\bf  X } & = \sum_{r=1}^R {\bf a}_r^{(1)} \circ {\bf a}_r^{(2)} \circ \cdots \circ {\bf a}_r^{(N)} \\
			&= \sum_{r=1}^R \gamma_r ({\bf u}_r^{(1)} \circ {\bf u}_r^{(2)} \circ \cdots \circ {\bf u}_r^{(N)} )
 \end{aligned}
\end{equation}
where ${\bf u}_r^{(n)} :={\bf a}_r^{(n)}/||{\bf a}_r^{(n)}||_2$, and character value $\gamma_r := \prod_{n=1}^N ||{\bf a}_r^{(n)}||_2 , r = 1, \cdots R$. Denote ${\bf U}^{(n)} := [ {\bf u}_1^{(n)}, {\bf u}_1^{(n)}, \ldots , {\bf u}_R^{(n)}]$, ${\bf U} := ({\bf U}^{(1)}, {\bf U}^{(2)}, \ldots, {\bf U}^{(N)})$, $\boldsymbol{\gamma} := (\gamma_1, \gamma_2, \ldots, \gamma_R)$, ${\bf A} := ({\bf A}^{(1)}, {\bf A}^{(2)}, \ldots, {\bf A}^{(N)})$

For fixed $n$, since vectors ${\bf a}^{(n)}_r$ ,$n= 1, \cdots N$ are interchangeable, identical distributions,  they can be modeled as independent for each other, zero-mean with the density function, 
 \begin{equation}\label{mix_model}
p_n({\bf a}^{(n)}_r) = \beta_{n}exp(-1/2 {\bf a}^{(n)}_r  {\bf R}_n^{-1}  {\bf a}^{(n)}_r - \mu_n || {\bf a}^{(n)}_r||_1)
\end{equation}
where $ \beta_n$ is the normalization constant, and the density is the product of the density of Gaussian (mean $0$, covariance ${\bf R}_n$) $1/\sqrt{\det(2 \pi {\bf R}_n )} \exp(-1/2  {\bf a}^{(n)}_r {\bf R}_n^{-1}  {\bf a}^{(n)}_r )$ and the double exponential density  $(\mu_n/2)^{I_n} \exp (- \mu_n || {\bf a}^{(n)}_r||_1)$.
% Note that the prior distribution of $ {\bf a}^{(n)}_r $, we have the intend that $\gamma_1 \ge \gamma_2 \cdots \ge \gamma_R$. 

  In addition $ {\bf a}^{(n)}_r, $ for $n= 1, 2, \ldots, N$ are assumed mutually independent. 
And since scale ambiguity is inherently present in the PARAFAC model, vectors  $ {\bf a}^{(n)}_r, $ for $n= 1, 2, \ldots, N$ are set to have equal power, that is,
 \begin{equation}\label{equi_power}
 \begin{aligned}
&\theta := \text{Tr}({\bf{R}}_1) = \text{Tr}({\bf{R}}_2) = \cdots = \text{Tr}({\bf{R}}_N) \\
&\mu := \mu_1 = \mu_2 = \ldots = \mu_N 
 \end{aligned}
\end{equation}
And for the computation tractable and simplicity, we assume that ${\bf{R}}_N$ are diagonal. 
Under these assumptions, the negative logarithm of the posterior distribution (up to a constant) is 
%proportional to $\exp(- L(\underline{\bf X}))$, with 
 \begin{equation}
 \begin{aligned}
&L(\underline{\bf X}) 
% = \frac{1}{2 \sigma^2} || (\underline{\bf Z} - \underline{\bf X}) {\varoast} \underline{{\boldsymbol {\Delta}}}  ||_F^2   \\
%			&+ \sum_{r= 1}^R \sum_{n= 1}^N  \left [   \frac{1}{2}  (({\bf a}^{(n)}_r)^T {\bf R}_n^{-1}{\bf a}^{(n)}_r)  + \mu    ||{\bf a}^{(n)}_r||_1 \right ]\\
%				&
				 = \frac{1}{2 \sigma^2} || (\underline{\bf Z} - \sum_{r=1}^R {\bf a}_r^{(1)} \circ {\bf a}_r^{(2)} \circ \cdots \circ {\bf a}_r^{(N)}) {\varoast} \underline{{\boldsymbol {\Delta}}}  ||_F^2 \\
				& + \sum_{r= 1}^R \sum_{n= 1}^N  \left [   \frac{1}{2}  (({\bf a}^{(n)}_r)^T {\bf R}_n^{-1}{\bf a}^{(n)}_r)  + \mu    ||{\bf a}^{(n)}_r||_1 \right ]\\
 \end{aligned}
\end{equation}
Correspondingly, the MAP estimator is  $ \underline{\bf  X } = {\bf A}^{(1)} \circ {\bf A}^{(2)}  \circ \cdots \circ {\bf A}^{(N)}  $, where ${\bf A}$ is the solution of
%$\underline{\hat {\bf X}} :=  \hat {\bf A}^{(1)} \circ \hat {\bf A}^{(2)}  \circ \cdots \circ \hat {\bf A}^{(N)}$,  where $\hat{\bf A}$ is the solution of 
% \begin{equation}
% \begin{aligned}
%&\underset{{\bf A} }{ \min  }\; \frac{1}{2 \sigma^2} || (\underline{\bf Z} - \sum_{r=1}^R {\bf a}_r^{(1)} \circ {\bf a}_r^{(2)} \circ \cdots \circ {\bf a}_r^{(N)}) {\varoast} \underline{{\boldsymbol {\Delta}}}  ||_F^2 \\
%				& + \sum_{r= 1}^R \sum_{n= 1}^N  \left [    \frac{1}{2}  (({\bf a}^{(n)}_r)^T {\bf R}_n^{-1}{\bf a}^{(n)}_r)  + \mu    ||{\bf a}^{(n)}_r||_1 \right ]\\
% \end{aligned}
%\end{equation}
%
%
%We can get a equivalent description,
 \begin{equation}\label{tensor_imputation_with_elastic_net}
 \begin{aligned}
&\underset{{\bf A} }{ \min  }\; \frac{1}{2 } || (\underline{\bf Z} - \sum_{r=1}^R {\bf a}_r^{(1)} \circ {\bf a}_r^{(2)} \circ \cdots \circ {\bf a}_r^{(N)}) {\varoast} \underline{{\boldsymbol {\Delta}}}  ||_F^2 \\
				&  + \sum_{r= 1}^R \sum_{n= 1}^N   \lambda\left [     \frac{1-\alpha}{2} (({\bf a}^{(n)}_r)^T {\bf R}_n^{-1}{\bf a}^{(n)}_r)  + \alpha  ||{\bf a}^{(n)}_r||_1 \right ]\\
 \end{aligned}
\end{equation}
where $\lambda = \sigma^2 (1+ \mu)$, $\alpha =\frac{\mu}{1+ \mu}$.
 Note that this regularization is elastic net as introduced in (\cite{elastic_net}). 
%(see in glmnet package)

% Note that this regularization is the elastic net in statistics  \cite{elastic_net}. 
%(see in glmnet package)

%%Columnwise 
%\section{Coordinate Descent Algorithm to Learn the Sparse and Low Rank PARAFAC Decomposition of the Observed Tensor with Missing Values}
\section{Block Coordinate Descent  for Optimization}
\label{sec:main_algorithm}
%Now we begin to derive the algorithm to solve (\ref{tensor_imputation_with_elastic_net}).
%We derive a novel coordinate descent algorithm to solve  (\ref{tensor_imputation_with_elastic_net}).
%%  the tensor imputation problem with the elastic net. 
% First, we separate the cost function in (\ref{tensor_imputation_with_elastic_net}) to the smooth part and the non-smooth part.
 In order to solve  (\ref{tensor_imputation_with_elastic_net}), we derive a coordinate descent-based algorithm, by firstly
%  the tensor imputation problem with the elastic net. 
separating the cost function in (\ref{tensor_imputation_with_elastic_net}) to the smooth part and the non-smooth part.
%Let ${\bf A} := ({\bf A}^{(1)}, {\bf A}^{(1)}, \ldots,  {\bf A}^{(N)})$, and let 
%For simplicity,  we define the same symbol function as in the article \cite{bcd}.
\begin{equation}\label{obj:bcd:f}
\begin{aligned}
f({\bf A}) := &\frac{1}{2 } || (\underline{\bf Z} - \sum_{r=1}^R {\bf a}_r^{(1)} \circ {\bf a}_r^{(2)} \circ \cdots \circ {\bf a}_r^{(N)}) {\varoast} \underline{{\boldsymbol {\Delta}}}  ||_F^2 \\
				&  + \sum_{r= 1}^R \sum_{n= 1}^N   \lambda\left [     \frac{1-\alpha}{2} (({\bf a}^{(n)}_r)^T {\bf R}_n^{-1}{\bf a}^{(n)}_r)  \right ]\\
 \end{aligned}
\end{equation}
%\begin{equation}\label{obj:bcd:r_n}
%\begin{aligned}
%r_n({\bf A}_n) :=  \lambda \alpha    \sum_{r= 1}^R ||{\bf a}^{(n)}_r||_1 \\ 
%\end{aligned}
%\end{equation}
\begin{equation}\label{obj:bcd:r}
\begin{aligned}
r({\bf A})&:= \sum_{n = 1}^N r_n({\bf A}_n)  =  \lambda \alpha    \sum_{r= 1}^R \sum_{n= 1}^N ||{\bf a}^{(n)}_r||_1 
\end{aligned}
\end{equation}
\begin{equation}\label{obj:bcd:F}
\begin{aligned}
F({\bf A}) := f( {\bf A}) + r({\bf A})
\end{aligned}
\end{equation}

The optimization problem in equation (\ref{tensor_imputation_with_elastic_net}) is solved iteratively by updating one part at a time with all other parts fixed. In detail, we cyclically minimize the columns ${\bf a}^{(n)}_r$  for $r=1, \ldots, R$ and $n=1, \ldots N$. 
For example, if we consider do minimization  about ${\bf a}_1^{(n)}$
  and fixed all the other columns of ${\bf A}$,
the following subproblem is then obtained:
\begin{equation}\label{tensor_imputation__elastic_net_for_column}
 \begin{aligned}
 F({\bf a}^{(n)}_1) : = & \frac{1}{2 } || (\underline{\bf Z} - \sum_{r=1}^R {\bf a}_r^{(1)} \circ {\bf a}_r^{(2)} \circ \cdots \circ {\bf a}_r^{(N)}) {\varoast} \underline{{\boldsymbol {\Delta}}}  ||_F^2 \\
				&  +  \lambda \left [     \frac{1-\alpha}{2} (({\bf a}^{(n)}_1)^T {\bf R}_n^{-1}{\bf a}^{(n)}_1)  + \alpha  ||{\bf a}^{(n)}_1||_1 \right ]\\
 \end{aligned}
\end{equation}
Setting $\underline{\bf W} =\underline{\bf Z} - \sum_{r=2}^R {\bf a}_r^{(1)} \circ {\bf a}_r^{(2)} \circ \cdots \circ {\bf a}_r^{(N)}$
\begin{equation}\label{tensor_imputation_elastic_net_for_column}
 \begin{aligned}
 F({\bf a}^{(n)}_1) : = & \frac{1}{2 } || (\underline{\bf W} -  {\bf a}_1^{(1)} \circ {\bf a}_1^{(2)} \circ \cdots \circ {\bf a}_1^{(N)}) {\varoast} \underline{{\boldsymbol {\Delta}}}  ||_F^2 \\
				&  +  \lambda \left [     \frac{1-\alpha}{2} (({\bf a}^{(n)}_1)^T {\bf R}_n^{-1}{\bf a}^{(n)}_1)  + \alpha  ||{\bf a}^{(n)}_1||_1 \right ]\\
 \end{aligned}
\end{equation}
%For this subproblem, we can alternative do minimization on one and fixed the another, e.g, do minimization on ${\bf a}^{(n)}_r$ and fixed   ${\bf a}^{(j)}_r, j \in [N]/\{n\}$ (where $[N] := \{1, 2, \ldots , N\}$), as known value. 
%For simplicity, we consider do minimization only on ${\bf a}^{(n)}_1$ and fixed ${\bf A}^{(m)}, m \in [N]/\{n\}$ and the other columns of ${\bf A}^{(n)}$. So the cost in  (\ref{tensor_imputation_with_elastic_net})  can be reduced to 
To make the calculation easier, we unfold the tensor to the matrix form. Let $ {\bf W}_{(n)} \in \mathbb R^{I_n \times I_1I_2\cdots I_{n-1} I_{n+1} \cdots I_N}$ denote the matrix of unfolding the tensor $\underline{\bf W}$ along its  mode $n$. And using the fact 
%if $\underline{\bf  X } = {\bf A}^{(1)} \circ {\bf A}^{(2)}  \circ \cdots \circ {\bf A}^{(N)}  $, then ${\bf X}_{(n)} = {\bf A}_n( {\bf A}_N \varodot \cdots \varodot  {\bf A}_{n+1} \varodot  {\bf A}_{n-1} \varodot \cdots \varodot  {\bf A}_{1})^T$.
$({\bf a}_1^{(1)} \circ {\bf a}_1^{(2)} \circ \cdots \circ {\bf a}_1^{(N)}) _{(n)}  = {\bf a}_1^{(n)} {\bf h}^T$, where  $ {\bf h} = {\bf a}_r^{\varoast _{-n}} :=  {\bf a}^{(N)}_r \varotimes \cdots \varotimes  {\bf a}^{(n+1)}_r \varotimes  {\bf a}^{(n-1)}_r \varotimes \cdots \varotimes  {\bf a}^{(1)}_r$. (\ref{tensor_imputation_elastic_net_for_column}) becomes
%
%
%To make a easy calculation , we  unfold the tensor to the matrix.  Let $ {\bf W}_{(n)} \in \mathbb R^{I_n \times I_1I_2\cdots I_{n-1} I_{n+1} \cdots I_N}$ denote the matrix of unfolding the tensor $\underline{\bf W}$ along its  mode $n$. And using the fact 
%%if $\underline{\bf  X } = {\bf A}^{(1)} \circ {\bf A}^{(2)}  \circ \cdots \circ {\bf A}^{(N)}  $, then ${\bf X}_{(n)} = {\bf A}_n( {\bf A}_N \varodot \cdots \varodot  {\bf A}_{n+1} \varodot  {\bf A}_{n-1} \varodot \cdots \varodot  {\bf A}_{1})^T$.
%$({\bf a}_1^{(1)} \circ {\bf a}_1^{(2)} \circ \cdots \circ {\bf a}_1^{(N)}) _{(n)}  = {\bf a}_1^{(n)} {\bf h}^T$, where  $ {\bf h} = {\bf a}_r^{\varoast _{-n}} :=  {\bf a}^{(N)}_r \varotimes \cdots \varotimes  {\bf a}^{(n+1)}_r \varotimes  {\bf a}^{(n-1)}_r \varotimes \cdots \varotimes  {\bf a}^{(1)}_r$. (\ref{tensor_imputation_elastic_net_for_column}) becomes
%Set ${\bf h} ={\bf a}_r^{\varoast _{-n}}  $. 
%Then we get the following subproblem about ${\bf a}_r^{(n)}$.
%can rewrite (\ref{tensor_imputation_elastic_net_for_column}) as 
\begin{equation}\label{tensor_imputation_elastic_net_for_column_unfold}
 \begin{aligned}
 F({\bf a}^{(n)}_1) : = & \frac{1}{2 } || ({\bf W}_{(n)} -  {\bf a}_1^{(n)}{\bf h}^T) {\varoast} {\boldsymbol {\Delta}}_{(n)}  ||_F^2 \\
				&  +  \lambda \left [     \frac{1-\alpha}{2} (({\bf a}^{(n)}_1)^T {\bf R}_n^{-1}{\bf a}^{(n)}_1)  + \alpha  ||{\bf a}^{(n)}_1||_1 \right ] \\
 \end{aligned}
\end{equation}
Denote ${\bf x}:= {\bf a}^{(n)}_1$ , ${\bf T}:= {\bf R}_n^{-1}$and drop out all subscript and superscript, we can rewrite (\ref{tensor_imputation_elastic_net_for_column_unfold}) as 
\begin{equation}\label{tensor_imputation_elastic_net_for_column_unfold_clean}
 \begin{aligned}
 F({\bf x}) : = & \frac{1}{2 } || ({\bf W} -  {\bf x}{\bf h}^T) {\varoast} {\boldsymbol {\Delta}} ||_F^2 \\
				&  +  \lambda \left [     \frac{1-\alpha}{2}{\bf x}^T {\bf T} {\bf x}  + \alpha  ||{\bf x}||_1 \right ] \\
 \end{aligned}
\end{equation}
which can be decomposed as 
\begin{equation}\label{tensor_imputation_elastic_net_fo_column_unfold_vector}
 \begin{aligned}
 F({\bf x}) & =  \sum_{i_n=1}^{I_n} [\frac{1}{2 } || {\boldsymbol {\delta}}_{i_n} \varoast {\bf w}_{i_n} -({\boldsymbol {\delta}}_{i_n} \varoast {\bf h} ){ x}_{i_n}||_2^2] \\
 				&  +  \lambda   \left [ \frac{1-\alpha}{2} {\bf x}^T {\bf T}{\bf x} + \alpha ||{\bf x}||_1 \right ]
 \end{aligned}
\end{equation}
where ${\bf w}_{i_n}^T$, ${\boldsymbol {\delta}}_{i_n}^T$, represent the ${i_n}$-th row of matrices $ {\bf W}$, ${\boldsymbol {\Delta}}$, respectively. 
%And ${\bf x} = [x_1, x_2, \ldots, x_{I_n}]^T$. 
%we can used the  coordinate descent method to solve (\ref{tensor_imputation_elastic_net_fo_column_unfold_vector}). 
Note that $F$ is strongly convex at ${\bf x}$, the optimal solution  
\begin{equation}\label{sol}
 \begin{aligned}
{\bf x}^*  = \arg \min_{\bf x} F({\bf x})
 \end{aligned}
\end{equation}
satisfies  the first order condition 
\begin{equation}\label{bcd:sol_condition}
 \begin{aligned}
{\bf 0} \in  \partial F({\bf x}^*)
 \end{aligned}
\end{equation}
The subgradient of $F$ at ${\bf x}$ is 
\begin{equation}\label{subgradient:F}
 \begin{aligned}
\partial F({\bf x}) &= ({\bf H}+\lambda  (1-\alpha) {\bf T}) {\bf x} - {\bf u } +\lambda   {\alpha}   \text{Sign}({\bf x})
 \end{aligned}
\end{equation}
where
\begin{equation}\label{H}       %开始数学环境
{\bf H} = \left[             %左括号
  \begin{array}{cccc}   %该矩阵一共3列，每一列都居中放置
    ||{\bf h}\varoast {\boldsymbol {\delta}}_1||_2^2& \quad &  \quad   \\  %第一行元素
     \quad 	& ||{\bf h} \varoast {\boldsymbol {\delta}}_2||_2^2&  \quad & \quad   \\  %第二行元素
    		 \quad  & \quad  & \ddots& \quad  \\  %第二行元素
 			 \quad &  \quad & \quad  &  ||{\bf h} \varoast {\boldsymbol {\delta}}_{I_n}||_2^2\\  %第二行元素
  \end{array}
\right     ]            %右括号
\end{equation}
%\begin{equation}\label{u} 
%{\bf u}= \left [({\bf h} \varoast {\boldsymbol {\delta}}_1)^T ({\bf w}_1\varoast {\boldsymbol {\delta}}_1), ({\bf h} \varoast {\boldsymbol {\Delta}}_2)^T ({\bf w}_2\varoast {\boldsymbol {\delta}}_2), \ldots ,({\bf h} \varoast {\boldsymbol {\delta}}_1)^T ({\bf w}_{I_n}\varoast {\boldsymbol {\delta}}_1) \right ]^T
%\end{equation}
\begin{equation}\label{u} 
{\bf u}= \left [({\bf h} \varoast {\boldsymbol {\delta}}_1)^T ({\bf w}_1\varoast {\boldsymbol {\delta}}_1), \ldots ,({\bf h} \varoast {\boldsymbol {\delta}}_{I_n})^T ({\bf w}_{I_n}\varoast {\boldsymbol {\delta}}_{I_n}) \right ]^T
\end{equation}
 and 
\begin{equation}\label{Sign}
  \text{Sign}(x)_i= \left \{
  \begin{aligned}
& 1 &\text{ if } x_i>0\\
& [-1,1] &\text{ if } x_i=0\\
& -1 &\text{ if } x_i<0\\
 \end{aligned}
 \right.
\end{equation}

%Now we decompose ${\bf T}:= {\bf T}^{(1)}+{\bf T}^{(2)}$, where ${\bf T}^{(1)} = \text{ diag }[T_{11}, T_{22}, \ldots,  T_{MM}]$, ${\bf T}^{(2)} := {\bf T} - {\bf T}^{(1)}$.
%, i.e. ${\bf T}^{(1)}$ is the diagonal matrix of ${\bf T}$ and ${\bf T}^{(2)}$ is the off-diagonal matrix of ${\bf T}$. 
Let 
\begin{equation}\label{d}
{\bf d} := \text{diag}( {\bf  H} +\lambda  (1-\alpha) {\bf T}):=(d_1, d_2, \ldots , d_{I_n})^T
\end{equation}
%\begin{equation}\label{subgradient:f_rewrite}
% \begin{aligned}
%\partial F({\bf x} ) =& ( {\bf  H} +\lambda  (1-\alpha) {\bf T}^{(1)})  {\bf x} +\lambda  (1-\alpha) {\bf T}^{(2)} {\bf x} \\
%& -{\bf u}  +\lambda   {\alpha} \text{Sign}({\bf x} )
% \end{aligned}
%\end{equation}
%Now 
%With some simple calculations, the optimal condition (\ref{bcd:sol_condition}) is equivalent to 
The optimal condition (\ref{bcd:sol_condition}) is actually equivalent to:
%\begin{equation}\label{bcd:sol_condition2}
% \begin{aligned}
%{\bf 0} \in \;&{\bf x} -  ( {\bf  H} +\lambda  (1-\alpha) {\bf T}^{(1)})^{-1} ({\bf u} -\lambda  (1-\alpha) {\bf T}^{(2)} {\bf x} )\\
%& + \lambda   {\alpha}( {\bf  H} +\lambda  (1-\alpha) {\bf T}^{(1)})^{-1} \text{Sign}({\bf x} )
% \end{aligned}
%\end{equation}
%For simplicity, 
%we can rewrite (\ref{bcd:sol_condition2})
\begin{equation}\label{bcd:sol_condition3}
 \begin{aligned}
 0 \in {\bf x} - {\bf u} \varoslash {\bf d} + \lambda   {\alpha}\text{Sign}({\bf x} )\varoslash {\bf d}
 \end{aligned}
\end{equation}
So the solution is
\begin{equation}\label{bcd:sol_condition3}
 \begin{aligned}
{\bf x} =  \mathcal T_{\lambda \alpha }({\bf u}) \varoslash {\bf d}
 \end{aligned}
\end{equation}
where $\mathcal T_v(t):= \text{sign}(t) \max\{|t|-v, 0\}$ is the soft thresholding operator.

%Now we can give the  algorithm \ref{alg:lrti_coordinate_descent}   to
%% learn the  sparse and low rank PARAFAC decomposition via the elastic net  problem 
% solve (\ref{tensor_imputation_with_elastic_net}).
We formally give the  algorithm \ref{alg:lrti_coordinate_descent} to
% learn the  sparse and low rank PARAFAC decomposition via the elastic net  problem 
 solve (\ref{tensor_imputation_with_elastic_net}).

\begin{algorithm}[tb]
%   \caption{Learning the  sparse and low rank PARAFAC decomposition via the elastic net regularization by coordinate decent}  
   \caption{Block coordinate descent method for solving (\ref{tensor_imputation_with_elastic_net})} 
   \label{alg:lrti_coordinate_descent}   
\begin{algorithmic}
   \STATE {\bfseries Input:}  Giving initial estimate  ${\bf A} := ({ \bf A}^{(1)}, { \bf A}^{(2)}, \ldots, { \bf A}^{(N)}  )$ of the factor matrix. 
   \REPEAT
  	 \STATE $\underline{\bf U }= \underline{\bf Z} - \underline{\bf X}$.
  	 \FOR{$r = 1$ {\bfseries to} $R$}
	 	\STATE update the factor pair ${(  { \bf a}^{(1)}_r,   { \bf a}^{(2)}_r,  \ldots, { \bf a}^{(N)}_r)}$
		\STATE 	$\underline{\bf W} =\underline{\bf U }+     { \bf a}^{(1)}_r \circ   { \bf a}^{(2)}_r \circ \cdots \circ  { \bf a}^{(N)}_r $
		\FOR{$n = 1$ {\bfseries to} $N$}
			\STATE 	Unfold ${\underline{\boldsymbol {{\boldsymbol {\Delta}}}}}$ and  $\underline{\bf W}$   	on the  mode $n$ into $ {\boldsymbol {{\boldsymbol {\Delta}}}}$ and  ${\bf W}$
			 \STATE 	Let ${\bf h} =     {\bf a}^{(N)}_r \varotimes \cdots \varotimes  {\bf a}^{(n+1)}_r \varotimes  {\bf a}^{(n-1)}_r \varotimes \cdots \varotimes  {\bf a}^{(1)}_r$, and  calculate  ${\bf u}$ and ${\bf d}$ as in equation (\ref{u}) and (\ref{d}), respectively.
%			 \STATE ${\bf d}:= [({\boldsymbol {\Delta}} \varoast ({\bf 1}_{I_n} {\bf h}^T)) \varoast ({\boldsymbol{\Delta}}  \varoast ({\bf 1}_{I_n} {\bf h}^T))] {\bf 1}_{I_1 \cdots I_{n-1} I_{n+1}\cdots I_N} +\lambda (1-\alpha) \text{diag}({\bf R}_n^{-1})$, \\
%			\STATE ${\bf u}= [({\boldsymbol {\Delta}} ({\bf 1}_{I_n} {\bf h}^T)) \varoast ( {\bf W} \varoast {\boldsymbol {\Delta}})] {\bf 1}_{I_1 \cdots I_{n-1} I_{n+1}\cdots I_N}$
			\STATE Update ${\bf a}^{(n)}_r$: 
				$ {\bf a}^{(n)}_r = \mathcal T_{\lambda \alpha }({\bf u}) \varoslash {\bf d}$
% Set ${\bf T}= {\bf R}_A^{-1}$, and $T_{m,m} =0$, for $m = 1, 2, \ldots, I_n$.
%			\FOR{$i_n = 1$ {\bfseries to} $I_n$}
%				 \STATE  ${\bf a}_r(i_n) =\mathcal T_{\frac{\lambda \alpha}{ d_{i_n}}}([u_m -\lambda (1- \alpha) ( {\bf T} {\bf a}_r)_{i_n}]/ d_{i_n})$
%			\ENDFOR
		\ENDFOR
		\STATE $\underline{\bf U } =\underline{\bf W} -     { \bf a}^{(1)}_r \circ   { \bf a}^{(2)}_r \circ \cdots \circ  { \bf a}^{(N)}_r $
   	\ENDFOR
	\STATE Let $\underline{\bf  X } = {\bf A}^{(1)} \circ {\bf A}^{(2)}  \circ \cdots \circ {\bf A}^{(N)} $.
   \UNTIL{Convergence}
\end{algorithmic}
\end{algorithm}

\section{Theory Analysis}
\label{sec:theory}
\subsection{Convergence results}
Since our algorithm is a instance of the unified algorithm 1 in (\cite{bcd}) , the convergence of algorithm \ref{alg:lrti_coordinate_descent} follows from Theory 2.8 and Theory 2.9 in (\cite{bcd}) and the proof is given in the Appendix section \ref{supp:theorem}.
We has  the following convergence theorem 
\begin{theorem}\label{th:convergence}
Let $\{ {\bf A}^k\} $ be the sequence generated by Algorithm  \ref{alg:lrti_coordinate_descent}, where  $\{ {\bf A}^k\} $  is the solution $({\bf A}^{(1)}, {\bf A}^{(2)}, \ldots , {\bf A}^{(N)} )$ in  the k-th iteration in the repeat loop. Assume that $\{ {\bf A}^k\} $ is bound. Then $\{ {\bf A}^k\} $ converges to  a critical point $\bar {\bf A}$, and the asymptotic convergence rates in Theory 2.9 in (\cite{bcd}) apply.
\end{theorem}

\subsection{Property of Elastic Regularization}
\subsubsection{$l_2$ regularization trend to find the  true rank}
First, the l2 regularization has the property to find the  true rank of tensor $\underline{\bf X}$ (\cite{rankreg}). To see this, note that when $\alpha = 0$( This can be achieved $\mu = 0$), the problem (\ref{tensor_imputation_with_elastic_net}) becomes 
 \begin{equation}\label{tensor_imputation_with_l2_norm}
 \begin{aligned}
\underset{{\bf A}, \; \underline{\bf X}}{\min} &\;\frac{1}{2 } || (\underline{\bf Z} - \underline{\bf X}) {\varoast} \underline{{\boldsymbol {\Delta}}}  ||_F^2 \\
				&  + \sum_{r= 1}^R \sum_{n= 1}^N   \lambda [     \frac{1}{2} ({\bf a}^{(n)}_r)^T {\bf R}_n^{-1}{\bf a}^{(n)}_r ]\\
				\text{s.t.}&\;	 \underline{\bf X} =  \sum_{r=1}^R {\bf a}_r^{(1)} \circ {\bf a}_r^{(2)} \circ \cdots \circ {\bf a}_r^{(N)}
 \end{aligned}
\end{equation}
And this problem is equivalent to the following problem by the following Proposition \ref{pro:l2equi}(its proof is provided in  the Appendix section \ref{supp:prop})
 \begin{equation}\label{tensor_imputation_with_l2_norm_standard_normalize}
 \begin{aligned}
&\underset{\tilde {\bf U}^{(n)}, \; \tilde \gamma_r ,\; \underline{\bf X}}{ \min} \frac{1}{2 } || (\underline{\bf Z} - \underline{\bf X}) {\varoast} \underline{{\boldsymbol {\Delta}}}  ||_F^2  + \sum_{r= 1}^R  \frac{\lambda N }{2}       ( {\tilde {\gamma}_r} )^{\frac{2}{N} } \\
				& \text{s.t. } \underline{\bf X} =\sum_{r=1}^R  \tilde \gamma_r ({\bf R}_1^{1/2} \tilde {\bf u}_r^{(1)} )  \circ \cdots \circ ( {\bf R}_N^{1/2} \tilde {\bf u}_r^{(N)})
\end{aligned}
\end{equation}
\begin{proposition}\label{pro:l2equi}
The solution of (\ref{tensor_imputation_with_l2_norm}) and (\ref{tensor_imputation_with_l2_norm_standard_normalize}) coincide, i.e. the optimal factors related by  
$ {\bf a}_r^{(N)} = (\tilde {\gamma}_r)^{\frac{1}{N}} {\bf R}_n^{1/2}\tilde {\bf u}_r^{(n)} $
\end{proposition}
%, \; n= 1, 2,\ldots ,N ,\; r = 1, 2, \ldots , R
To further stress the capability of (\ref{tensor_imputation_with_l2_norm}) to produce a low-rank approximate  tensor, consider transform  (\ref{tensor_imputation_with_l2_norm_standard_normalize}) once more by rewritten it in the constrained-error form 
 \begin{equation}\label{tensor_imputation_with_l2_norm_standard_5}
 \begin{aligned}
&\underset{\tilde {\bf U}^{(n)} , \;\tilde {\boldsymbol {\gamma} } ,\; \underline{\bf X}}  { \min} 
				|| \tilde {\boldsymbol{ \gamma}} ||_{\frac{2}{N} }   \\
				& \text{ s.t. }     {\underline {\bf X}} = \sum_{r=1}^R  \tilde \gamma_r ({\bf R}_1^{1/2} \tilde {\bf u}_r^{(1)} )  \circ \cdots \circ ( {\bf R}_N^{1/2} \tilde {\bf u}_r^{(N)}) \\
				&|| \left (\underline{\bf Z} -  \underline {\bf X} \right ) {\varoast} \underline{{\boldsymbol {\Delta}}}  ||_F^2 \le \eta
\end{aligned}
\end{equation}
where $ \tilde {\boldsymbol {\gamma} } := ( \tilde \gamma_1,  \tilde \gamma_2, \ldots,  \tilde \gamma_R)$, $|| \tilde {\boldsymbol{ \gamma}} ||_{\frac{2}{N} }   := (\sum_{r=1}^R  |\tilde \gamma_r |^{\frac{2}{N}})^{\frac{N}{2}}$. For any vlue of $\eta$ there exists a corresponding Lagrange multiplier $\lambda$ such that (\ref
%{tensor_imputation_with_l2_norm_standard_3}
{tensor_imputation_with_l2_norm_standard_normalize})  and (\ref{tensor_imputation_with_l2_norm_standard_5}) yield the same solution. And the $l_{\frac{N}{2}}$-norm $|| \tilde {\boldsymbol{ \gamma}} ||_{\frac{2}{N} } $ in (\ref{tensor_imputation_with_l2_norm_standard_5}) produces a sparse vector $\tilde {\boldsymbol{ \gamma}}$ when minimized (\cite{lp-sparse-solution}). Note that the sparsity in vector $\tilde {\boldsymbol{ \gamma}}$ implies the low rank of ${\underline {\bf X}}$.

The above arguments imply that when properly choose the parameter $\lambda$, the $l_2$ regularization problem (\ref{tensor_imputation_with_l2_norm}) can find the true rank of the tensor $\underline{\bf X}$. i.e. when we give a overestimated rank $R> R^*$(where $R^*$ is the true rank of  of tensor $\underline{\bf X}$) of tensor $\underline{\bf X}$, the perfect solution of  (\ref{tensor_imputation_with_l2_norm}) will shrink the redundant $R- R^*$ columns of the factor matrix ${\bf A}$ to $0$.

\subsubsection{$l_1$ regularization trend to find the true sparse structure of the factor matrix}
Note that when there is only $l_1$ regularization
%( which can be achieved  when the factor comes from
% \begin{equation}
%p_n({\bf a}^{(n)}_r) = \beta_{n}exp( - \mu || {\bf a}^{(n)}_r||_1)
%\end{equation}
%, then $\lambda = \mu$, $\alpha = 0$), the problem (\ref{tensor_imputation_with_elastic_net}) becomes 
% \begin{equation}\label{tensor_imputation_with_elastic_net_l1_norm}
% \begin{aligned}
%\underset{{\bf A} }{\min} &\; \frac{1}{2 } || (\underline{\bf Z} - \sum_{r=1}^R {\bf a}_r^{(1)} \circ {\bf a}_r^{(2)} \circ \cdots \circ {\bf a}_r^{(N)}) {\varoast} \underline{{\boldsymbol {\Delta}}}  ||_F^2 \\
%				&  + \sum_{r= 1}^R \sum_{n= 1}^N   \lambda  ||{\bf a}^{(n)}_r||_1\\
% \end{aligned}
%\end{equation}
%Which can be written in a compact form
 \begin{equation}\label{tensor_imputation_with_elastic_net_l1_norm_compact}
 \begin{aligned}
\underset{{\bf A} }{\min} &\; \frac{1}{2 } || (\underline{\bf Z} -  {\bf A}^{(1)} \circ {\bf A}^{(2)} \circ \cdots \circ {\bf A}^{(N)}) {\varoast} \underline{{\boldsymbol {\Delta}}}  ||_F^2 \\
				&  +  \sum_{n= 1}^N   \lambda  ||\text{vec}({\bf A}^{(n)})||_1\\
 \end{aligned}
\end{equation}
For each mode factor matrix ${\bf A}^{(n)}$,  it is a standard lasso problem, which implies that the solution $({\bf A}^{(n)})^*$ is sparse. Note that if  when the true factor matrix $({\bf A}^{(n)})^*$ is sparse, with the properly choosed $\lambda$, we can reveal the true  sparsity structure in  $({\bf A}^{(n)})^*$. For application, the  sparsity structure of  $({\bf A}^{(n)})^*$ can help us to make meaningful explanation on its mode $n$, which standards for some attributes of the considered problem.

\subsubsection{The elastic net give us a flexible model to find the true data structure in tensor}
	Combine the $l_1$ regularization and $l_2$ regularization we get the elastic net regularization, which can combine the advantages of both $l_1$ regularization and $l_2$ regularization, i.e. to find the true (low) rank  and closed to the true sparse factor matrix of PARAFAC decomposition of an observed tensor data.  It is helpful to understand the structure of the data and reveal the faces of the objects. 

\subsubsection{Estimate of the covariance matrix}	
To run the algorithm \ref{alg:lrti_coordinate_descent}, the covariance must be postulated as a priori, or replace by their sample estimates. And often we don't know the priori, so  the proper sample estimates are very important, since it provides reasonable scaling in each dimension such that the algorithm performs well. Similarly in section C Covariance estimation in article (\cite{rankreg}), we can  bridge  the covariance matrix with its kernel counterparts.  Since that when the model is give in (\ref{mix_model})  
%\begin{equation}
%p_n({\bf a}^{(n)}_r) = \beta_{n}exp(-1/2 {\bf a}^{(n)}_r {\bf R}_n^{-1}  {\bf a}^{(n)}_r - \mu_n || {\bf a}^{(n)}_r||_1)
%\end{equation} 
the covariance of the the factor ${\bf a}^{(n)}_r$ is hard to evaluate, do not has an analytical expression. We can jump out the obstacles just assume our model is Gaussian now,
\begin{equation}\label{gaussian_model}
p_n({\bf a}^{(n)}_r) = \beta_{n}exp(-1/2 {\bf a}^{(n)}_r {\bf R}_n^{-1}  {\bf a}^{(n)}_r )
\end{equation} 

Define the kernel similarity matrix in mode $n$ as
\begin{equation}\label{kernel_similarity_matrix}
{\bf K}_n(i,j) := \mathbb E \; {\bf X}_{(n)} (i,:) ({\bf X}_{(n)} (j,:))^T
\end{equation} 
i.e. the expectation of the inner product of the $i$-th slice and $j$-th slice of $ {\bf X}$ in mode $n$.
With some calculation(see in the supplementary material), we can get 
\begin{equation}\label{kernel_similarity_matrix_sol2}
{\bf K}_n = R\theta^{N-1} {\bf R}_n
\end{equation} 
and 
\begin{equation}\label{X_F_norm2}
\mathbb E ||\underline{\bf X}||_F^2 = \text{Tr}({\bf K}_n)=  R\theta^{N}
\end{equation}
From this we can get the covariance matrix estimate(just drop out the expectation)
\begin{equation}\label{cov_estimate}
 \begin{aligned}
&\theta = (\frac{ ||\underline{\bf X}||_F^2 }{R})^{\frac{1}{N}}\\
&{\bf K}_n(i,j) = {\bf X}_{(n)} (i,:) ({\bf X}_{(n)} (j,:))^T\\
&{\bf R}_n =\frac{ {\bf K}_n}{R \theta^{N-1}}
 \end{aligned}
\end{equation}

\section{Solution Path by Warm Start in Practical Application}
\label{sec:solution_path}
%stochastic implemention in the algorithm.
%
%extrapolation to accelate
%
%Solution path by warm start. 
%, \; n= 1, 2,\ldots ,N 
In this section, we will give a more practical algorithm which  help us how to chose  the proper sparse and low rank PARAFAC approximation of the observed tensor with missing values. From the Lagrange theory, the problem (\ref{tensor_imputation_with_elastic_net}) is equivalent to 
 \begin{equation}\label{tensor_imputation_with_elastic_net_constrained}
 \begin{aligned}
 &\underset{{\bf A} }{\min} \;\frac{1}{2 } || (\underline{\bf Z} - \sum_{r=1}^R {\bf a}_r^{(1)} \circ {\bf a}_r^{(2)} \circ \cdots \circ {\bf a}_r^{(N)}) {\varoast} \underline{{\boldsymbol {\Delta}}}  ||_F^2 \\
				& \text{s.t.  } \sum_{r= 1}^R \sum_{n= 1}^N  \left [     \frac{1-\alpha}{2} (({\bf a}^{(n)}_r)^T {\bf R}_n^{-1}{\bf a}^{(n)}_r)  + \alpha  ||{\bf a}^{(n)}_r||_1 \right ] \le \eta\\
 \end{aligned}
\end{equation}

For any value of $\lambda$ there exists a corresponding Lagrange multiplier $\eta$ such that (\ref{tensor_imputation_with_elastic_net})   and (\ref{tensor_imputation_with_elastic_net_constrained}) yield the same solution. Although there is no explicit equation to relate $\lambda$ and $\eta$, we know that small $\lambda$ implies  large $\eta$, large ball-like solution space; large $\lambda$ implies  small $\eta$, small ball-like solution space. So we can fixed the relation $\alpha$  of $l_1$ and $l_2$ regularization. And start from a small $\lambda$ toward to large $\lambda$, which can shrink the dense solution to the sparse solution(when $\lambda$ is very large, the solution is $0$). And this forms a solution path, the solutions are close to each other when $\lambda$ are close to each other, this means when we start from the previous solution, the algorithm can find quickly the close $\lambda$ solution. And for the lucky $\alpha$, we can meet the true data structure, i.e. the right sparse structure in the factor matrix. From our numerical test, the situation is just as we imagined. However, on the time when we find the right  data structure, i.e. the right sparse structure in the factor matrix, the recovery error is not so low, since when $\lambda$ is large, the regularization will shrink the solution close to origin, and the best recovery error is achieved when  $\lambda$  is small. This phenomenon inspired us a efficient method: using algorithm \ref{alg:lrti_coordinate_descent}  to generate a solution path, and for a special solution in the solution path, fixed it sparse structure,  solve a constrained problem with small $\lambda$ to get small recovery error --- for fixed $\alpha$, we use Algorithm \ref{alg:lrti_coordinate_descent} to calculate a solution path, i.e. the solution generated by  Algorithm \ref{alg:lrti_coordinate_descent}
from a sequence increasing $\lambda$ (e.g $(10^{-10}, 10^{-9}, \ldots, 10^{10})$. 

For a special solution ${\bf A}:=  ({ \bf A}^{(1)}, { \bf A}^{(2)}, \ldots, { \bf A}^{(N)}  )$
, where ${ \bf A}^{(n)} \in \mathbb R^{I_n \times R}$,  calculate ${\bf U}$ and $\boldsymbol{\gamma}$ as in (\ref{cp_decompositon}), reorder $\boldsymbol{\gamma}$ such that $\gamma_1\ge \gamma_2 \ge \gamma_{R_1} >\gamma_{R_1+1}= \cdots =\gamma_{R} =0$ , reorder the columns of matrix ${ \bf A}^{(n)}$ and  ${ \bf U}^{(n)}$ such that they coincide with the $\boldsymbol{\gamma}$. Let ${\bf B}:=  ({ \bf B}^{(1)}, { \bf B}^{(2)}, \ldots, { \bf B}^{(N)}  )$ , where ${ \bf B}^{(n)} \in \mathbb R^{I_n \times R_1}$ is sub-matrix of ${ \bf A}^{(n)}$, which is formed by its first $R_1$ columns. Definde the sparse structure matrix ${\bf S}:=   ({ \bf S}^{(1)}, { \bf S}^{(2)}, \ldots, { \bf S}^{(N)}  )$ , where ${ \bf S}^{(n)} \in \mathbb R^{I_n \times R_1}$ of  the factor matrix ${\bf B}$ as 
 \begin{equation}\label{sparse-structure-matrix}
{ \bf S}^{(n)}(i_n,r) := \left \{
 \begin{aligned}
  0 & \text{ if } |{\bf B}^{(n)}(i_n, r)| \le \epsilon \\ 
  1 & \text{ if } |{\bf B}^{(n)}(i_n, r)| > \epsilon
 \end{aligned}
 \right .
\end{equation}
where $i \in  [I_n]$, $j \in [R_1]$ and $\epsilon$ is a prespecified small number( e.g. $\epsilon = 10^{-9}$). 

For notation simplicity, we still use ${\bf A}$ denote ${\bf B}$ and $R$ denote $R_1$. If the true sparse structure matrix is  ${\bf S}$, to solve (\ref{tensor_imputation_with_elastic_net}),  it is equivalent to solve the constrained problem, 
 \begin{equation}\label{tensor_imputation_with_elastic_net_sparse_constrained}
 \begin{aligned}
 \underset{{\bf A}}{\min} &\;\frac{1}{2 } || (\underline{\bf Z} - \sum_{r=1}^R {\bf a}_r^{(1)} \circ {\bf a}_r^{(2)} \circ \cdots \circ {\bf a}_r^{(N)}) {\varoast} \underline{{\boldsymbol {\Delta}}}  ||_F^2 \\
				&+ \sum_{r= 1}^R \sum_{n= 1}^N \lambda \left [     \frac{1-\alpha}{2} (({\bf a}^{(n)}_r)^T {\bf R}_n^{-1}{\bf a}^{(n)}_r)  + \alpha  ||{\bf a}^{(n)}_r||_1 \right ]  \\
				\text{s.t. } &{ \bf A}^{(n)}(i_n,r) = 0, \text{ if } { \bf S}^{(n)}(i_n,r) =0  \\
				&\text{ for } i_n \in [I_n], \; r \in [R]
 \end{aligned}
\end{equation}
This problem can solved by algorithm \ref{alg:lrti_coordinate_descent} just renew ${ \bf A}^{(n)}(i_n,r)$ if  ${ \bf S}^{(n)}(i_n,r) =1$, and set all ${ \bf A}^{(n)}(i_n,r)=0$ if  ${ \bf S}^{(n)}(i_n,r) =0$. We formally give the Algorithm \ref{alg:lrti_coordinate_descent_sparse}.

\begin{algorithm}[tb]
%   \caption{Low-rank tensor imputation with elastic net regularization by coordinate decent given the sparse structure matrix} 
     \caption{Sparse constrained coordinate descent method for
 solving (\ref{tensor_imputation_with_elastic_net_sparse_constrained})}
   \label{alg:lrti_coordinate_descent_sparse}    
\begin{algorithmic}
   \STATE  Giving initial estimate  ${\bf A} := ({ \bf A}^{(1)}, { \bf A}^{(2)}, \ldots, { \bf A}^{(N)}  )$ of the factor matrix and the sparse structure matrix ${\bf S}:=   ({ \bf S}^{(1)}, { \bf S}^{(2)}, \ldots, { \bf S}^{(N)}  )$. All the steps are same as in Algorithm \ref{alg:lrti_coordinate_descent} except we replace the inner most for loop by the following sparse constrained form.
			\FOR{$i_n = 1$ {\bfseries to} $I_n$}
			 	\IF{${ \bf S}^{(n)}(i_n,r) =1$}
%				 	 \STATE  ${\bf a}_r(i_n) =\mathcal T_{\frac{\lambda \alpha}{ d_{i_n}}}(u_m)/ d_{i_n})$
					  \STATE  ${\bf a}_r(i_n) =\mathcal T_{\lambda \alpha}(u_m)/ d_{i_n})$

				\ENDIF
			\ENDFOR
\end{algorithmic}
\end{algorithm}

And the convergence theory still holds since we only update a subset elements of factor matrix ${\bf A}$.

To conclude, we give our finally algorithm \ref{alg:solution_path} aim to find the true data structure and low recovery error.  For convenience, we call the solution path method using  algorithm \ref{alg:lrti_coordinate_descent} SPML, and the correspond using  algorithm \ref{alg:lrti_coordinate_descent_sparse} SPMS where the solution comes from SPML. 

\begin{algorithm}[!tb]   
\caption{Solution path method}   
\label{alg:solution_path}   
\begin{algorithmic}[1]  
%\State Let $ \underline{\bf Z}(m,n,p)= 0$, if ${\boldsymbol {\Delta}}(m,n,p) =0$.
\STATE  Given initial estimate  ${\bf A}_0 := ({ \bf A}^{(1)}, { \bf A}^{(2)}, \ldots, { \bf A}^{(N)}  )$ of the factor matrix. Given $\alpha$,   a increasing $\lambda$ sequence $ (\lambda_1, \lambda_2, \ldots, \lambda_L)$, and an small  $\lambda_s$ and $\epsilon$. 
\FOR{$l=1$ {\bfseries to} $L$}
	\STATE  Using Algorithm \ref{alg:lrti_coordinate_descent}  with $\lambda_l$ and $\alpha$ and  ${\bf A}_{l-1}$ as initial estimate of factor matrix to generate solution ${\bf A}_l$. 
	\STATE  Let ${\bf A} :={\bf A}_l$, calculate ${\bf U}$ and $\boldsymbol{\gamma}$ as in (\ref{cp_decompositon}), reorder $\boldsymbol{\gamma}$ such that $\gamma_1\ge \gamma_2 \ge \gamma_{R_1} >\gamma_{R_1+1}= \cdots =\gamma_{R} =0$ , reorder the columns of matrix ${ \bf A}^{(n)}$ and  ${ \bf U}^{(n)}$ such that they coincide with the $\boldsymbol{\gamma}$. Let ${\bf B}:=  ({ \bf B}^{(1)}, { \bf B}^{(2)}, \ldots, { \bf B}^{(N)}  )$ , where ${ \bf B}^{(n)} \in \mathbb R^{I_n \times R_1}$ is sub-matrix of ${ \bf A}^{(n)}$, which is formed by its first $R_1$ columns. Definde the sparse structure matrix ${\bf S}:=   ({ \bf S}^{(1)}, { \bf S}^{(2)}, \ldots, { \bf S}^{(N)}  )$ , where ${ \bf S}^{(n)} \in \mathbb R^{I_n \times R_1}$ of  the factor matrix ${\bf B}$ as 
 \begin{equation}\label{Sparsity matrix}
{ \bf S}^{(n)}(i_n,r) := \left \{
 \begin{aligned}
  0 & \text{ if } |{\bf B}^{(n)}(i_n, r)| \le \epsilon \\ 
  1 & \text{ if } |{\bf B}^{(n)}(i_n, r)| > \epsilon
 \end{aligned}
 \right .
\end{equation}
where $i \in  [I_n]$, $j \in [R_1]$. Let $s$ denote the number of $0$ in ${\bf S}$,
 \begin{equation}\label{Sparsity number}
 \begin{aligned}
  s = \sum_{n=1}^N   \text{ number of } 0 \text{ in } {\bf S}^{(n)}
 \end{aligned}
\end{equation}
	\STATE  Denote the sparse structure matrix  ${\bf S}_l := {\bf S}$, and the rank and number of $0$ in ${\bf S}$ pair $(R_l, s_l) := (R_1, s)$, where  $R_1$ and $s$ comes from step 4. Let ${\bf A}_l = {\bf B}$
	\STATE  When the sparse structure changed, i.e. ${\bf S}_l$ is differrent with ${\bf S}_{l-1}$, using Algorithm \ref{alg:lrti_coordinate_descent_sparse}  with $\lambda_s$ , $\alpha$ , ${\bf A}_l$ as initial estimate of factor matrix and  ${\bf S}_l $ as the sparse data matrix to generate solution ${\bf B}_l$
	\STATE  Choose the best factor matrix ${\bf A}$ from the solution path $\{{\bf B}_l \}$ with low recovery error and the most sparse structure.
\ENDFOR
\end{algorithmic}  
\end{algorithm}

\section{Generalize the algorithm to large scale problem by using stochastic block coordinate descent algorithm by using admax method}
\label{sec:main_algorithm_random}
For the real application, 
%such as in deep learning, one want to represent  a convolution kernel by the PARAFAC model to reduce the parameters and save the computation cost. In such a situation,
the data often have  a  large data dimnsion.
the above algorithm may not compute fast, since for each iteration, we have $O(I_1*I_2*\cdots *I_N)$ cost of operations. Fortunately, we can fixed this problem to some extent by using a stochastic scheme.

 Note that for renew ${\bf a}^{(n)}_r$, we use all information to renew it, it it not necessary. To gain the insights. Now consider the problem ${\bf x} = a{\bf b}$, where ${\bf x} \in \mathbb R^{M \times 1}$ ,  $a \in R$ and ${\bf b} \in \mathbb R^{M \times 1}$. Now suppose we know the true $ {\bf b}$, and we need to update $a$, the form the above algorithm, we renew $a$ by the formula(Now we drop out the regularity and suppose their is no noise and no missing entries) $ a = \frac{\langle {\bf x}, {\bf b} \rangle}{\langle {\bf b}, {\bf b}\rangle}$. Note that we can also update $a$ by random choose a subset ${\bf s} \subseteq \{ 1, 2, \ldots, M\}$, and update $ a = \frac{\langle {\bf x_{\bf s}}, {\bf b_{\bf s}} \rangle}{\langle {\bf b_{\bf s}}, {\bf b_{\bf s}}\rangle}$. While when there are noise, such that we observed ${\bf z} = {\bf x} + {\bf e}$, $ \tilde a = \frac{\langle {\bf z }, {\bf b} \rangle}{\langle {\bf b}, {\bf b}\rangle} =\frac{\langle {\bf x} +  {\bf e} , {\bf b} \rangle}{\langle {\bf b}, {\bf b}\rangle}  = a +  \frac{\langle  {\bf e} , {\bf b} \rangle}{\langle {\bf b}, {\bf b}\rangle}$. Using a subset update $ \hat a = \frac{\langle {\bf z }_{\bf s}, {\bf b}_{\bf s} \rangle}{\langle {\bf b}_{\bf s}, {\bf b}_{\bf s}\rangle} =\frac{\langle {\bf x}_{\bf s} +  {\bf e}_{\bf s} , {\bf b}_{\bf s} \rangle}{\langle {\bf b}_{\bf s}, {\bf b}_{\bf s} \rangle}  = a +  \frac{\langle  {\bf e}_{\bf s} , {\bf b}_{\bf s} \rangle}{\langle {\bf b}_{\bf s}, {\bf b}_{\bf s}\rangle}$. Suppose that $e \sim \mathcal N(0 ,\sigma^2 {\bf I})$ is a white noise. then $\mathbb E \tilde a = \mathbb E \hat a $, and $\text{Var} (\tilde a)  = \sigma^2/||{\bf b}||_2^2$, $\text{Var} (\hat a)  = \sigma^2/||{\bf b}_{\bf s}||_2^2$. From this simple case, it inspire us that we can use a subset of ${\bf b}$ to update $a$, because it will lose some accuracy, higher variance, we should it a stochastic update scheme to renew $a$. Its key idea is keep the memory of the older information of the right direction and also take notice of the current of new information. 

Now we begin to deduce our stochastic  block coordinate descent method by using the Adamax scheme. Follow the same roads in the Section \ref{sec:main_algorithm}.

For example, if we consider do minimization  about ${\bf a}_1^{(n)}$, fixed all the other columns of ${\bf A}$, and we choose subset ${\bf s}_k \subseteq  \{ 1, 2 , \ldots, I_k\}$ for each mode $k = 1,2, \ldots ,N$. 
Denote 
$$  \underline{\bf Z}_n :=  \underline{\bf Z}({\bf s}_1, \ldots, {\bf s}_{n-1},:, {\bf s}_{n+1}, \ldots, {\bf s}_{N})$$
the following subproblem is then obtained:
\begin{equation}\label{tensor_imputation__elastic_net_for_column_subset}
 \begin{aligned}
 &F({\bf a}^{(n)}_1) : = \\
 & \frac{1}{2 } || (\underline{\bf Z}_n - \sum_{r=1}^R 
 {\bf a}_r^{(1)}({\bf s}_1) \circ \cdots \circ {\bf a}_{r}^{(n-1)}({\bf s}_{n-1})  \circ {\bf a}_r^{(n)} \\
 &\circ {\bf a}_r^{(n+1)}({\bf s}_{n+1}) \circ \cdots \circ {\bf a}_r^{(N)}({\bf s}_{N}) ) 
 {\varoast} \underline{{\boldsymbol {\Delta}}}_n ||_F^2 \\
				&  +  \lambda \left [     \frac{1-\alpha}{2} (({\bf a}^{(n)}_1)^T {\bf R}_n^{-1}{\bf a}^{(n)}_1)  + \alpha  ||{\bf a}^{(n)}_1||_1 \right ]\\
 \end{aligned}
\end{equation}
Setting $\underline{\bf W}_n =\underline{\bf Z}_n - \sum_{r=2}^R {\bf a}_r^{(1)}({\bf s}_1) \circ \cdots \circ {\bf a}_{r}^{(n-1)}({\bf s}_{n-1})  \circ {\bf a}_r^{(n)} \circ {\bf a}_r^{(n+1)}({\bf s}_{n+1}) \circ \cdots \circ {\bf a}_r^{(N)}({\bf s}_{N}) $
\begin{equation}\label{tensor_imputation_elastic_net_for_column_subset}
 \begin{aligned}
 F({\bf a}^{(n)}_1) : = & \frac{1}{2 } || (\underline{\bf W}_n -  
 {\bf a}_1^{(1)}({\bf s}_1) \circ \cdots \circ {\bf a}_{1}^{(n-1)}({\bf s}_{n-1})  \circ {\bf a}_1^{(n)} \\
 &\circ {\bf a}_1^{(n+1)}({\bf s}_{n+1}) \circ \cdots \circ {\bf a}_1^{(N)}({\bf s}_{N}) 
 ) {\varoast} \underline{{\boldsymbol {\Delta}}}_n  ||_F^2 \\
				&  +  \lambda \left [     \frac{1-\alpha}{2} (({\bf a}^{(n)}_1)^T {\bf R}_n^{-1}{\bf a}^{(n)}_1)  + \alpha  ||{\bf a}^{(n)}_1||_1 \right ]\\
 \end{aligned}
\end{equation}
%For this subproblem, we can alternative do minimization on one and fixed the another, e.g, do minimization on ${\bf a}^{(n)}_r$ and fixed   ${\bf a}^{(j)}_r, j \in [N]/\{n\}$ (where $[N] := \{1, 2, \ldots , N\}$), as known value. 
%For simplicity, we consider do minimization only on ${\bf a}^{(n)}_1$ and fixed ${\bf A}^{(m)}, m \in [N]/\{n\}$ and the other columns of ${\bf A}^{(n)}$. So the cost in  (\ref{tensor_imputation_with_elastic_net})  can be reduced to 
To make the calculation easier, we unfold the tensor to the matrix form. Let $ {\bf W}_{(sn)} \in \mathbb R^{I_n \times |{\bf s_1}|  |{\bf s_2}|\cdots  |{\bf s_{n-1}}|  |{\bf s_{n+1}}| \cdots  |{\bf s_N}|}$ denote the matrix of unfolding the tensor $\underline{\bf W}_n$ along its  mode $n$. And using the fact 
%if $\underline{\bf  X } = {\bf A}^{(1)} \circ {\bf A}^{(2)}  \circ \cdots \circ {\bf A}^{(N)}  $, then ${\bf X}_{(n)} = {\bf A}_n( {\bf A}_N \varodot \cdots \varodot  {\bf A}_{n+1} \varodot  {\bf A}_{n-1} \varodot \cdots \varodot  {\bf A}_{1})^T$.
$ {\bf a}_1^{(1)}({\bf s}_1) \circ \cdots \circ {\bf a}_{1}^{(n-1)}({\bf s}_{n-1})  \circ {\bf a}_1^{(n)} 
 \circ {\bf a}_1^{(n+1)}({\bf s}_{n+1}) \circ \cdots \circ {\bf a}_1^{(N)}({\bf s}_{N} )
 = {\bf a}_1^{(n)} {\bf h_s}^T$, where  $ {\bf h}_s :=  {\bf a}^{(N)}_1 ({\bf s}_{N}) \varotimes \cdots \varotimes  {\bf a}^{(n+1)}_1({\bf s}_{n+1}) \varotimes  {\bf a}^{(n-1)}_1({\bf s}_{n-1}) \varotimes \cdots \varotimes  {\bf a}^{(1)}_1({\bf s}_{1})$. (\ref{tensor_imputation_elastic_net_for_column_subset}) becomes
\begin{equation}\label{tensor_imputation_elastic_net_for_column_unfold_subset}
 \begin{aligned}
 F({\bf a}^{(n)}_1) : = & \frac{1}{2 } || ({\bf W}_{(sn)} -  {\bf a}_1^{(n)}{\bf h}_s^T) {\varoast} {\boldsymbol {\Delta}}_{(sn)}  ||_F^2 \\
				&  +  \lambda \left [     \frac{1-\alpha}{2} (({\bf a}^{(n)}_1)^T {\bf R}_n^{-1}{\bf a}^{(n)}_1)  + \alpha  ||{\bf a}^{(n)}_1||_1 \right ] \\
 \end{aligned}
\end{equation}
Denote ${\bf x}_s:= {\bf a}^{(n)}_1$ , ${\bf T}:= {\bf R}_n^{-1}$and drop out all subscript and superscript, and to denote it is update by  the subset information, we add a subscript $s$ for each symbol,   we can rewrite (\ref{tensor_imputation_elastic_net_for_column_unfold_subset}) as 
\begin{equation}\label{tensor_imputation_elastic_net_for_column_unfold_clean_subset}
 \begin{aligned}
 F({\bf x}_s) : = & \frac{1}{2 } || ({\bf W}_s -  {\bf x}_s{\bf h}_s^T) {\varoast} {\boldsymbol {\Delta}}_s ||_F^2 \\
				&  +  \lambda \left [     \frac{1-\alpha}{2} {\bf x}_s^T {\bf T} {\bf x}_s  + \alpha  ||{\bf x}_s||_1 \right ] \\
 \end{aligned}
\end{equation}
which can be decomposed as 
\begin{equation}\label{tensor_imputation_elastic_net_fo_column_unfold_vector_subset}
 \begin{aligned}
 F({\bf x}_s) & =  \sum_{i_n=1}^{I_n} [\frac{1}{2 } || {\boldsymbol {\delta}}_{s,i_n} \varoast {\bf w}_{s, i_n} -({\boldsymbol {\delta}}_{s, i_n} \varoast {\bf h} ){ x}_{s, i_n}||_2^2] \\
 				&  +  \lambda   \left [ \frac{1-\alpha}{2} {\bf x}_s^T {\bf T}{\bf x}_s + \alpha ||{\bf x}_s||_1 \right ]
 \end{aligned}
\end{equation}
where ${\bf w}_{s, i_n}^T$, ${\boldsymbol {\delta}}_{s, i_n}^T$, represent the ${i_n}$-th row of matrices $ {\bf W}_s$, ${\boldsymbol {\Delta}}_s$, respectively. 
Note that $F$ is strongly convex at ${\bf x}_s$, the optimal solution  
\begin{equation}\label{sol}
 \begin{aligned}
{\bf x}_s^*  = \arg \min_{{\bf x}_s} F({\bf x}_s)
 \end{aligned}
\end{equation}
satisfies  the first order condition 
\begin{equation}\label{bcd:sol_condition_subset}
 \begin{aligned}
{\bf 0} \in  \partial F({\bf x}_s^*)
 \end{aligned}
\end{equation}
The subgradient of $F$ at ${\bf x}$ is 
\begin{equation}\label{subgradient:F_subset}
 \begin{aligned}
\partial F({\bf x}_s) &= ({\bf H}_s+\lambda  (1-\alpha) {\bf T}) {\bf x}_s - {\bf u }_s +\lambda   {\alpha}   \text{Sign}({\bf x}_s)
 \end{aligned}
\end{equation}
where
\begin{equation}\label{H_subset}       %开始数学环境
{\bf H}_s = \left[             %左括号
  \begin{array}{cccc}   %该矩阵一共3列，每一列都居中放置
    ||{\bf h}_s\varoast {\boldsymbol {\delta}}_{s,1}||_2^2& \quad &  \quad   \\  %第一行元素
     \quad 	& ||{\bf h}_s \varoast {\boldsymbol {\delta}}_{s,2}||_2^2&  \quad & \quad   \\  %第二行元素
    		 \quad  & \quad  & \ddots& \quad  \\  %第二行元素
 			 \quad &  \quad & \quad  &  ||{\bf h}_s \varoast {\boldsymbol {\delta}}_{s, I_n}||_2^2\\  %第二行元素
  \end{array}
\right     ]            %右括号
\end{equation}
\begin{equation}\label{u_subset} 
\begin{aligned}
&{\bf u}_s= [({\bf h}_s \varoast {\boldsymbol {\delta}}_{s,1})^T ({\bf w}_{s,1}\varoast {\boldsymbol {\delta}}_{s,1}), \ldots ,\\
&({\bf h}_s \varoast {\boldsymbol {\delta}}_{s,I_n})^T ({\bf w}_{s,I_n}\varoast {\boldsymbol {\delta}}_{s, I_n})]^T
\end{aligned}
\end{equation}
 and 
\begin{equation}\label{Sign_subset}
  \text{Sign}(x)_i= \left \{
  \begin{aligned}
& 1 &\text{ if } x_i>0\\
& [-1,1] &\text{ if } x_i=0\\
& -1 &\text{ if } x_i<0\\
 \end{aligned}
 \right.
\end{equation}

Let 
\begin{equation}\label{d_subset}
{\bf d}_s := \text{diag}( {\bf  H}_s +\lambda (1-\alpha) {\bf T}):=(d_{s,1}, d_{s,2}, \ldots , d_{s,I_n})^T
\end{equation} 
The optimal condition (\ref{bcd:sol_condition_subset}) is actually equivalent to:
\begin{equation}\label{bcd:sol_condition3_subset}
 \begin{aligned}
 0 \in {\bf x}_s - {\bf u}_s \varoslash {\bf d}_s + \lambda   {\alpha}\text{Sign}({\bf x}_s )\varoslash {\bf d}_s
 \end{aligned}
\end{equation}
So the solution is
\begin{equation}\label{bcd:sol_condition3_subset}
 \begin{aligned}
{\bf x}_s =  \mathcal T_{\lambda \alpha }({\bf u}_s) \varoslash {\bf d}_s
 \end{aligned}
\end{equation}
where $\mathcal T_v(t):= \text{sign}(t) \max\{|t|-v, 0\}$ is the soft thresholding operator.

Note that although we solve for ${\bf a}_1^{(n)}$ using a subset information, however, we want it has the approximate the same effect as the original solution. Note that  from equation (\ref{bcd:sol_condition3}) and (\ref{bcd:sol_condition3_subset}), ${\bf u}_s \varoslash {\bf d}_s \approx {\bf u} \varoslash {\bf d}$ is no problem, however ${\lambda \alpha } \varoslash {\bf d}_s \not \approx {\lambda \alpha }  \varoslash {\bf d}$. So the correct update should be like 
\begin{equation}\label{bcd:sol_condition3_subset_correct}
 \begin{aligned}
{\bf x}_s =  \mathcal T_{\lambda \alpha  \varoslash  {\bf d}}({\bf u}_s \varoslash {\bf d}_s)
 \end{aligned}
\end{equation}
Since we do not have ${\bf d}$ in the subset update form, we now deduce an approximate formula for equation (\ref{bcd:sol_condition3_subset_correct}). 
Note that when $ \underline{{\boldsymbol {\Delta}}}= \underline{\bf 1}$, ${\bf d} = ||{\bf h}||_2^2{ \bf 1 }+ \lambda (1- \alpha)\alpha \text{diag}({\bf T})$  and ${\bf d}_s = ||{\bf h}_s||_2^2{ \bf 1 }+ \lambda (1- \alpha)\alpha \text{diag}({\bf T})$, so ${\bf d}  \approx \text{mean}( \frac{ ||{\bf h}||_2^2{ \bf 1 }+ \lambda (1- \alpha)\alpha \text{diag}({\bf T})}{{\bf d}_s}){\bf d}_s  $. 
Then the update can be approximated by 
\begin{equation}\label{bcd:sol_condition3_subset_correct_approx}
 \begin{aligned}
{\bf x}_s =  \mathcal T_{\lambda \alpha \text{ mean}( \frac{{\bf d}_s}{ ||{\bf h}||_2^2{ \bf 1 }+ \lambda (1- \alpha)\alpha \text{diag}({\bf T})})}({\bf u}_s )\varoslash {\bf d}_s
 \end{aligned}
\end{equation}
Where the $ ||{\bf h}||_2^2 = \prod_{i=1, i\ne n}^N ||{\bf a}_1^{(i)}||_2^2$.  

Note that for the Admax algorithm(\cite{adam}), it calculate the gradient ${\bf g}_1^{(n)} $ of $F$ at ${\bf a}_1^{(n)(t-1)}$, then begin to update, but now we have the true solution of $\min_{ {\bf a}_1^{(n)}} F ({\bf a}_1^{(n)}$, we just set ${\bf g}_1^{(n)}  = {\bf a}_1^{(n)(t-1)} -{\bf x}_s $

We formally give the algorithm \ref{alg:lrti_coordinate_descent_subset} to
% learn the  sparse and low rank PARAFAC decomposition via the elastic net  problem 
 solve (\ref{tensor_imputation_with_elastic_net}). For convenience, we call the solution path method using  algorithm \ref{alg:lrti_coordinate_descent_subset} SPMLR, and the correspond using  algorithm \ref{alg:lrti_coordinate_descent_sparse} SPMSR where the solution comes from SPMLR. 

\begin{algorithm}[tb]
%   \caption{Learning the  sparse and low rank PARAFAC decomposition via the elastic net regularization by coordinate decent}  
   \caption{ Adamax Stochastic Block Coordinate descent method for solving (\ref{tensor_imputation_with_elastic_net})} 
   \label{alg:lrti_coordinate_descent_subset}   
\begin{algorithmic}
   \STATE {\bfseries Input:}  Giving initial estimate  ${\bf A} := ({ \bf A}^{(1)}, { \bf A}^{(2)}, \ldots, { \bf A}^{(N)}  )$ of the factor matrix.  And the batch size $(s_1, s_2, \ldots, s_N)$ for each mode, where $s_n \in \mathbb N^+$ and $ 1\le s_n \le I_n$.   ${\bf M} := ({ \bf M}^{(1)}, { \bf M}^{(2)}, \ldots, { \bf M}^{(N)}  )$ be the first moment estimate matrix, which has the same shape as ${\bf A}$ and is initialized to $0$. 
 Let  ${\bf U}_{alg} \in \mathbb R^{N \times R}$ be exponentially weighted infinity norm matrix   and it  is initialized to $0$. Let $h_{normsqure} \in \mathbb R^{R}$, where $h_{normsqure}(r) = \prod_{n=1}^N ||{\bf a}_r^{(n)}||_2^2$. Given the exponential decay rates $\beta_1, \; \beta_2 \in [0,1)$(default, $\beta_1=0.9, \; \beta_2=0.9999$), and given step-size $\alpha_{alg}$(which is decreased to  $0.2*\alpha_{alg}$ once the relative recovery error $\frac{||\underline{\bf Z} - \underline{\bf X} ||_F}{||\underline{\bf Z}||_F}$ begin to increase). Let ${\underline {\bf U}}$ be the same shape as ${\underline {\bf Z}}$ and initialized to $0$.
  \STATE Set $t=0$. 
   \REPEAT
   	\FOR{$r =1$ {\bfseries to} $R$}
	 \STATE $h_{normsqure}(r) = \beta_1 h_{normsqure}(r) +(1-\beta_1)\prod_{n=1}^N ||{\bf a}_r^{(n)}||_2^2$. 
	\ENDFOR
%  	 \STATE $\underline{\bf U }= \underline{\bf Z} - \underline{\bf X}$.
	\STATE Randomly (uniformly) select a set with size $s_n$ ${\bf s}_n \in \{ 1, 2, \ldots, I_n\}$  for each mode. 
	\STATE Upadate $\underline{\bf U } = \underline{\bf Z} - \underline{\bf X} $ only for following need used parts. i.e. $\underline{\bf U } ({\bf s}_1, \ldots, {\bf s}_{n-1},: ,  {\bf s}_{n+1}, \ldots, {\bf s}_{N})$ for  $n = 1, \ldots, N$.
  	 \FOR{$r = 1$ {\bfseries to} $R$}
	 	\STATE update the factor pair ${(  { \bf a}^{(1)}_r,   { \bf a}^{(2)}_r,  \ldots, { \bf a}^{(N)}_r)}$
		\FOR{$n = 1$ {\bfseries to} $N$}
			\STATE Let ${\bf s} = ({\bf s}_1, \ldots,{\bf s}_{n-1},:, {\bf s}_{n}, \ldots, {\bf s}_N )$
			\STATE 	$\underline{\bf W}({\bf s})  =\underline{\bf U }({\bf s}) +     {\bf a}_r^{(1)}({\bf s}_1) \circ \cdots \circ {\bf a}_{r}^{(n-1)}({\bf s}_{n-1})  \circ {\bf a}_r^{(n)} \circ {\bf a}_r^{(n+1)}({\bf s}_{n+1}) \circ \cdots \circ {\bf a}_r^{(N)}({\bf s}_{N} )$
			\STATE 	Unfold ${\underline{\boldsymbol {{\boldsymbol {\Delta}}}}}({\bf s})  $ and  $\underline{\bf W}({\bf s}) $   	on the  mode $n$ into $ {\boldsymbol {{\boldsymbol {\Delta}}}}_s$ and  ${\bf W}_s$
			 \STATE 	Let $ {\bf h}_s :=  {\bf a}^{(N)}_r ({\bf s}_{N}) \varotimes \cdots \varotimes  {\bf a}^{(n+1)}_r({\bf s}_{n+1}) \varotimes  {\bf a}^{(n-1)}_r({\bf s}_{n-1}) \varotimes \cdots \varotimes  {\bf a}^{(1)}_r({\bf s}_{1})$, and  calculate  ${\bf u}_s$ and ${\bf d}_s$ as in equation (\ref{u_subset}) and (\ref{d_subset}), respectively.
			\STATE Calculate 
			$ \tau = \text{ mean}( {{\bf d}_s} \varoslash({ \frac{h_{normsquare}(r) }{ ||{\bf a}^{(n)}_r||_2^2}{ \bf 1 }+ \lambda (1- \alpha)\alpha \text{diag}({\bf T})}))$\\
			 $ \bar {\bf a}_r^{(n)}  =\mathcal T_{\lambda \alpha \tau}({\bf u}_s )\varoslash {\bf d}_s$
			 \STATE Set the gradients ${\bf g}_r^{(n)} = {\bf a}_r^{(n)} - \bar {\bf a}_r^{(n)} $
			 \STATE Update biased first moment estimate 
			 $ {\bf M}^{(n)}(:,r) = \beta_1  {\bf M}^{(n)}(:,r) +  (1- \beta_1 ) {\bf g}_r^{(n)}$
			 \STATE Update the exponentially weighted infinity norm
			 $ {\bf U}_{alg}(n ,r) =\max (\beta_2 {\bf U}_{alg}(n ,r), ||{\bf g}_r^{(n)}	||_2)$
			  \STATE Update $	 {\bf a}_r^{(n)} =  {\bf a}_r^{(n)} - (\alpha_{alg} / (1- \beta_1^t)) {\bf M}^{(n)}(:,r) /{\bf U}_{alg}(n ,r) $\\
			   $	 {\bf a}_r^{(n)}(\bar {\bf a}_r^{(n)} =0)=0$
			
% Set ${\bf T}= {\bf R}_A^{-1}$, and $T_{m,m} =0$, for $m = 1, 2, \ldots, I_n$.
%			\FOR{$i_n = 1$ {\bfseries to} $I_n$}
%				 \STATE  ${\bf a}_r(i_n) =\mathcal T_{\frac{\lambda \alpha}{ d_{i_n}}}([u_m -\lambda (1- \alpha) ( {\bf T} {\bf a}_r)_{i_n}]/ d_{i_n})$
%			\ENDFOR
		\ENDFOR
		 \STATE Upadate $\underline{\bf U } = \underline{\bf W} - {\bf a}_r^{(1)} \circ \cdots \circ {\bf a}_r^{(N)}$ only for following need used parts. i.e. $\underline{\bf U } ({\bf s}_1, \ldots, {\bf s}_{n-1},: ,  {\bf s}_{n+1}, \ldots, {\bf s}_{N})$ for  $n = 1, \ldots, N$.
%		 $\underline{\bf U }({\bf s}) =\underline{\bf W}({\bf s}) -   {\bf a}_r^{(1)}({\bf s}_1) \circ \cdots \circ {\bf a}_{r}^{(n-1)}({\bf s}_{n-1})  \circ {\bf a}_r^{(n)} \circ {\bf a}_r^{(n+1)}({\bf s}_{n+1}) \circ \cdots \circ {\bf a}_r^{(N)}({\bf s}_{N} )$
   	\ENDFOR
	\STATE $t =t+1$
   \UNTIL{Convergence}
\end{algorithmic}
\end{algorithm}

%\section{Convergence Results}

%In this section, we should prove 
The convergence theorem for Algorithm \ref{alg:lrti_coordinate_descent_subset} is not clear. It belong to the class of the multiconvex stochastic scheme, we make the following conjecture.
%but I don't know how to prove it. 
\begin{conjecture}
Let $\{ {\bf A}^k\} $ be the sequence generated by Algorithm  \ref{alg:lrti_coordinate_descent_subset}, where  $\{ {\bf A}^k\} $  is the solution $({\bf A}^{(1)}, {\bf A}^{(2)}, \ldots , {\bf A}^{(N)} )$ in  the k-th iteration in the repeat loop. Assume that $\{ {\bf A}^k\} $ is bound. Then $\{ {\bf A}^k\} $ converges to  a critical point $\bar {\bf A}$.
\end{conjecture}

%\twocolumn[
\section{Numerical Test}
\label{sec:numerical_test}

\subsection{Simulated Data}
In this work, we generate the synthetic dataset to test the performance of the proposed algorithm. Without loss of generality, the simulation is done on the 3-way tensor  data where $I := I_1=I_2=I_3$ using the Bayesian model as described in section 2. For simplicity, we assume that ${\bf R}_n$ is  a diagonal matrix, where its diagonal comes from uniform distribution. i.e.
%We now simulate the synthetic data to test the performance of our algorithm. For simplicity, we just do simulation on the 3-way tensor  data where $I := I_1=I_2=I_3$.
%We simulate our data using the Bayesian model in section 2.
%For simplicity, we test on that ${\bf R}_n$ is  a diagonal matrix, where its diagonal comes from uniform distribution. i.e.
\newpage
 \begin{equation}\label{Sparsity number}
 \begin{aligned}
 {\bf R}_1(1,1) \sim U(101,131) \\
 {\bf R}_1(i,i) \sim U(1,31), \; i = 2, \ldots, I_1 \\
 {\bf R}_2(1,1) \sim U(1001,1021) \\
 {\bf R}_2(i,i) \sim U(1,21), \; i = 2, \ldots, I_2 \\
 {\bf R}_3(1,1) \sim U(10001,10011) \\
 {\bf R}_3(i,i) \sim U(1,11), \; i = 2, \ldots, I_3 \\
 \end{aligned}
\end{equation}
Note that when the  ${\bf R}_n$ is diagonal, 
 \begin{equation}\label{mix_model_element}
  \begin{aligned}
p_n({\bf a}^{(n)}_r(i)) =& \beta_{n}exp(-1/2  \;{\bf a}^{(n)}_r(i) {\bf R}_n(i,i)^{-1}  {\bf a}^{(n)}_r(i )\\
&- \mu | {\bf a}^{(n)}_r(i)|)
 \end{aligned}
\end{equation}
We can sample from the distribution by the rejection sample method. 
To produce the sparse factors, we set a gate $g$. When sampling out the factor matrix ${\bf A}$, we make it sparse by letting all its elements whose absolute value less than gate $g$ to $0$.  

\subsubsection{Solution path}
\label{section:no noise and no missing}
We firstly use a simple example to investigate the behavior of the solution path, by setting $\mu = 0.1$ and $g=0.5$. $(I_1, I_2, I_3) :=(6 , 6, 6)$. $R=2$  and the true number of zeros in ${\bf A}$ is $12$.  And we test our algorithm with initial estimate rank $R_{\text{est}} = 3$ , $\alpha = 0.2$, and the covariance ${\bf R}$ is estimated from equation (\ref{cov_estimate}) (where we only estimate the diagonal of  ${\bf R}$, if one diagonal of ${\bf R}$ less than $1e-8$, we just set it to $1e-8$ to ensure stability). The initial estimate matrix ${\bf A}_0$  comes from standard Gaussian noise.  The numbers of  max iterations in Algorithm \ref{alg:lrti_coordinate_descent} and Algorithm \ref{alg:lrti_coordinate_descent_sparse} are $200$ and  $100$, respectively.

%\begin{table}[!tb]
\begin{table}[ht]
\caption{Solution path by Algorithm \ref{alg:solution_path} on the simultated  data when $\alpha = 0.2$. R: rank, NZS: the number of zeros in the current solution, NZT: the number of the  zeros such that the zero appears in the current solution and the true solution;
IS1: which algorithm is used?  1--algorithm \ref{alg:lrti_coordinate_descent} and 0-- algorithm \ref{alg:lrti_coordinate_descent_sparse}; 
}
\vskip 0.1in
\label{table:solution path}
\begin{center} 
\begin{small}
\begin{sc}
\begin{tabular}{ccccccc}
\hline
\abovespace\belowspace
%$\lambda $&		N2 & N3 &N4 &N5&N6&N7\\
$\lambda $&		R & NZS &NZT &IS1&rel\_err&iters\\
 \hline
 \abovespace
1e-10&     3&     3&     2& 1&8.8e-06& 38 \\  
1e-08&     3&     3&     2& 0&6.4e-06& 2 \\  
1e-04&     3&     3&     2& 1&1.8e-05& 56 \\  
1e-08&     3&     3&     2& 0&3.5e-06& 3 \\  
1e-01&     3&     8&     4& 1&3.7e-03& 200 \\  
1e-08&     3&     8&     4& 0&2.3e-05& 100 \\  
1e+00&     2&     5&     5& 1&2.3e-03& 200 \\  
1e-08&     2&     5&     5& 0&4.4e-06& 12 \\  
1e+01&     2&     7&     7& 1&1.3e-02& 200 \\  
1e-08&     2&     7&     7& 0&3.2e-06& 15 \\  
2e+01&     2&     8&     8& 1&2.5e-02& 200 \\  
1e-08&     2&     8&     8& 0&4.0e-06& 16 \\  
2e+02&     2&    10&    10& 1&1.9e-01& 108 \\  
1e-08&     2&    10&    10& 0&3.3e-06& 17 \\  
4e+02&     2&    11&    11& 1&3.2e-01& 63 \\  
1e-08&     2&    11&    11& 0&4.4e-06& 8 \\  
5e+02&     2&    12&    12& 1&3.8e-01& 53 \\  
1e-08&     2&    12&    12& 0&1.6e-06& 9 \\  
7e+02&     1&     3&     3& 1&6.1e-01& 56 \\  
1e-08&     1&     3&     3& 0&5.4e-01& 10 \\  
\belowspace
2e+03&     0&     0&     0& 1&1.0e+00& 6 \\ \hline
\end{tabular}
\end{sc}
\end{small}
\end{center}
\vskip -0.1in
\end{table}

In table \ref{table:solution path} (the full table can be found in the Appendix Section \ref{supp:solution_path}), we see that as $\lambda$ increases, the algorithm will shrink the dense solution to $0$, and decrease the rank $3$ to $0$. For all the the resulting factors, their sparsity increases and rank decreases along the solution path as more elements in the factors shrink to 0. When the true rank $2$ is found, the sparsity structure of the synthetic data is also found by the algorithm gradually, as indicated by the increasing number of zero elements found in the current solution and  in the true solution, as shown in column NZS and NZT. When $\lambda = 500$, the algorithm \ref{alg:lrti_coordinate_descent} finds the true sparse structure of the  true solution ${\bf A}$. However the relative recovery error using $\lambda = 500$ is quite high ($3.78e-01$). While Algorithm \ref{alg:lrti_coordinate_descent_sparse} can decrease the  relative recovery error to $1.62e-6$ (find the true solution ${\bf A}$) in just 9 iterations. 

We observed that when $\lambda$ is small, the algorithm converges very fast. But when $\lambda$ becomes large, the iterations will exceed the max number of iterations we set despite the fact that we used the warm start strategy. Yet it still contains the useful information for the true solution and makes it possible to decrease the  relative recovery error from $3.78e-01$ to $1.62e-6$ only in 9 iterations.

%\twocolumn[
\subsubsection{Model Performance with Noise and No Missing Values}
\label{subsubsection:random_test_behavior}
In this section, we test our algorithm with some baseline algorithms. The first algorithm is the CP-ALS method(which is provided by the Tensor Toolbox for MATLAB (\cite{tensortoolbox}) , and we reimplement it in python), and the second algorithm is the LRTI algorithm in (\cite{rankreg}) (where the parameter $\mu$ in LRTI algorithm chosen by $\mu = 0.01 *\mu_{\text{max}}$ ).

To compare the numerical results, we use standard score metric to measure how well the ground truth is recovered by a CP decomposition in those case where the true factors are known  (\cite{comparealg}). The score between two rank-one tensors $\underline{\bf X} = {\bf a} \circ {\bf b} \circ {\bf c}$ and  $\underline{\bf Y} = {\bf p} \circ {\bf q} \circ {\bf r}$ is defined as:
\begin{equation}
\text{score}( \underline{\bf X} , \underline{\bf Y}) = \frac{{\bf a} ^T{\bf p} }{||{\bf a}||\; || {\bf p}||} \times \frac{{\bf b} ^T{\bf q} }{||{\bf b}|| \; || {\bf q}||}   \times \frac{{\bf c} ^T{\bf r} }{||{\bf c}|| \; || {\bf r}||}
\end{equation}
For $R > 1$, we first sort the components such that $\gamma_1 \ge \gamma_2 \ge \cdots \ge \gamma_R$, then average the scores for all pairs of components. 

After the above simple synthetic data experiment, we further design a more complex synthetic scheme with random noise added. Let $(I_1, I_2, I_3)= (20,20,20)$, $R=10$, the number of test times is set to $50$, and the initial estimate of rank is $10$, and $\mu = 0.1$, $g=0.5$. 
We add noise to the data with the SNR of 20dB, where the SNR is defined by $\text{SNR} := 10\log_{10} \frac{\text{Var}(\underline{\bf X}) }{\text{Var}(\underline{\bf E})}$, and $\sigma^2 = \text{Var}(\underline{\bf E})$. And the is no missing values $\underline{\boldsymbol{\Delta}}= \underline{\bf 1}$, since the CP-ALS algorithm cannot not hander the missing value entries.

 For each test in our solution path method( SPML, SPMS. we don't test on SPMLR, SPMSR algorithms, since the dimension is only $20$. but  the solution are chosen by the same procedure in the next subsection), we select the solution in the solution according to the following steps: (1) Choose the estimate rank be 
 the smallest found rank which not less than the true rank; 
(2) For those solutions whose rank is  the estimate rank found by (1), choose the last solution which satisfies that the  zeros in the current solution are also the zeros in the true solution, and the number of zeros in the current solution is not greater than the true number of zeros, if no such solution, go to step (3);
(3) For those solutions whose rank is  the estimate rank found by (1), choose  the last solution which the number of zeros in the current solution less than the true number of zeros.
(4) If not such a solution, we can choose a best solution from the path by ourself.

And for each algorithm,  we use two ways to initialize the ${\bf A}$, the first one is  to use   the random standard Gaussian $\mathcal N(0,1)$, we call it random initialzation, and the second is to set ${\bf A}^{(n)}$ to be the leading $R$ left singular vectors of the mode-n unfolding, ${\bf Z}_{(n)} \varoast {\boldsymbol{\Delta}}_{(n)}$, and we call this nvecs initialization. Now we give our test result for both random initialization and  nvecs initialization.

%\twocolumn[

The results are summarized in the table \ref{table:with-noise-no-missing20_10_10in_paper}. 

%\twocolumn[
%\newpage
%]

\begin{table}[!tb]
\caption{50 times random test behavior when $(I_1, I_2, I_3)= (20,20,20)$, $R=10$  where the data with noise ( SNR = 20dB) and  without missing entries and each algorithm using random initialization (given in above) and  nvecs initialization(given below).
And SPML stands for the Solution Path Method using  algorithm \ref{alg:lrti_coordinate_descent} (Look algorithm), and And SPMS stands for the Solution Path Method using  algorithm \ref{alg:lrti_coordinate_descent_sparse} (Sparse constrained  Algorithm). And SPML(0.8) standes for that $\alpha = 0.8$. And SCORES(STD) stand for the mean score and the its standard deviation. 
REL\_ERR stands for the mean of the full data relative error   $\frac{||\underline{\bf Z} - \underline{\bf X} ||_F}{||\underline{\bf Z}||_F}$ subtract  the true relative error  $\frac{||\underline{\bf E}||_F}{||\underline{\bf Z}||_F}$.
TFNZ stands for the total number of zeros found by the selected solution;  
TFTNZ stands for  total number of the zeros  such that the zero
appears in both the selected solution and the true solution;
TNZ stands for  total number of the zeros  in the true solution.
TIME stands for the the mean cost cpu time in 50 random testes, note that for each $\alpha$, SPML and SPMS are computed in together, so the time are both the  SPML and SPMS in each test.
}
\label{table:with-noise-no-missing20_10_10in_paper}
\vskip 0.15in
\begin{center} 
%\begin{small}
\footnotesize
\begin{sc}
\begin{tabular}{ccccccccc}
\hline
\abovespace\belowspace
%$\lambda $&		N2 & N3 &N4 &N5&N6&N7\\
Method&		Scores(std) &rel\_err & TFNZ & TFTNZ & TNZ &$\frac{\text{TFNZ }}{\text{TNZ}}$ &$\frac{\text{TTNZ }}{\text{TNZ}}$   & time \\
 \hline
\abovespace\belowspace

random \quad &initialization \quad &\quad &\quad &\quad &\quad &\quad &\quad \\
\hline 
\abovespace
CP-AlS	&	0.566(0.360)	& 	-2.436e-03	&  	     0	&  	     0	&  	  5078	&  	0.000	&  	0.000	& 	0.6\\ 
LRTI	&	0.413(0.278)	& 	1.918e-02	&  	     0	&  	     0	&  	  5078	&  	0.000	&  	0.000	& 	1.2\\ 
SPML(0)	&	0.579(0.400)	& 	1.824e-02	&  	     0	&  	     0	&  	  5009	&  	0.000	&  	0.000	& 	73.7\\ 
SPMS(0)	&	0.648(0.387)	& 	1.435e-02	&  	     0	&  	     0	&  	  5009	&  	0.000	&  	0.000	& 	73.7\\ 
SPML(0.2)	&	0.922(0.141)	& 	1.178e-02	&  	  2225	&  	  2130	&  	  5009	&  	0.444	&  	0.425	& 	74.5\\ 
SPMS(0.2)	&	0.967(0.121)	& 	-3.351e-03	&  	  2225	&  	  2225	&  	  5009	&  	0.444	&  	0.444	& 	74.5\\ 
SPML(0.8)	&	0.934(0.164)	& 	2.771e-03	&  	  3003	&  	  2933	&  	  5009	&  	0.600	&  	0.586	& 	67.3\\ 
SPMS(0.8)	&	0.952(0.160)	& 	-2.646e-03	&  	  3003	&  	  2989	&  	  5009	&  	0.600	&  	0.597	& 	67.3\\ 
SPML(0.98)	&	0.946(0.144)	& 	1.956e-03	&  	  3168	&  	  3083	&  	  5009	&  	0.632	&  	0.615	& 	67.1\\ 
\belowspace
SPMS(0.98)	&	0.960(0.140)	& 	-3.324e-03	&  	  3168	&  	  3112	&  	  5009	&  	0.632	&  	0.621	& 	67.1\\ 
\hline
\abovespace\belowspace
nvecs \quad & initialization\quad &\quad &\quad &\quad &\quad &\quad &\quad \\
\hline 
CP-AlS	&	0.739(0.289)	& 	-9.452e-03	&  	     0	&  	     0	&  	  5078	&  	0.000	&  	0.000	& 	0.4\\ 
LRTI	&	0.647(0.240)	& 	1.632e-02	&  	     2	&  	     0	&  	  5078	&  	0.000	&  	0.000	& 	2.4\\ 
SPML(0)	&	0.784(0.324)	& 	7.721e-03	&  	     0	&  	     0	&  	  5009	&  	0.000	&  	0.000	& 	109.0\\ 
SPMS(0)	&	0.832(0.285)	& 	4.766e-03	&  	     0	&  	     0	&  	  5009	&  	0.000	&  	0.000	& 	109.0\\ 
SPML(0.2)	&	0.883(0.226)	& 	1.360e-02	&  	  2274	&  	  2146	&  	  5009	&  	0.454	&  	0.428	& 	113.0\\ 
SPMS(0.2)	&	0.927(0.224)	& 	-2.580e-03	&  	  2274	&  	  2241	&  	  5009	&  	0.454	&  	0.447	& 	113.0\\ 
SPML(0.8)	&	0.888(0.251)	& 	4.696e-03	&  	  3090	&  	  2905	&  	  5009	&  	0.617	&  	0.580	& 	105.3\\ 
SPMS(0.8)	&	0.909(0.252)	& 	-1.160e-03	&  	  3090	&  	  2967	&  	  5009	&  	0.617	&  	0.592	& 	105.3\\ 
SPML(0.98)	&	0.897(0.234)	& 	4.672e-03	&  	  3344	&  	  3037	&  	  5009	&  	0.668	&  	0.606	& 	104.4\\ 
\belowspace
SPMS(0.98)	&	0.921(0.230)	& 	-1.498e-03	&  	  3344	&  	  3084	&  	  5009	&  	0.668	&  	0.616	& 	104.4\\ 
\hline
\end{tabular}
\end{sc}
%\end{small}
\end{center}
\vskip -0.1in
\end{table}
%]
%   \twocolumn[
%\newpage
%]
% \twocolumn[

\begin{figure}[ht]
\vskip 0.2in
\begin{center}
\includegraphics[width=\columnwidth,height=0.6\textwidth ]{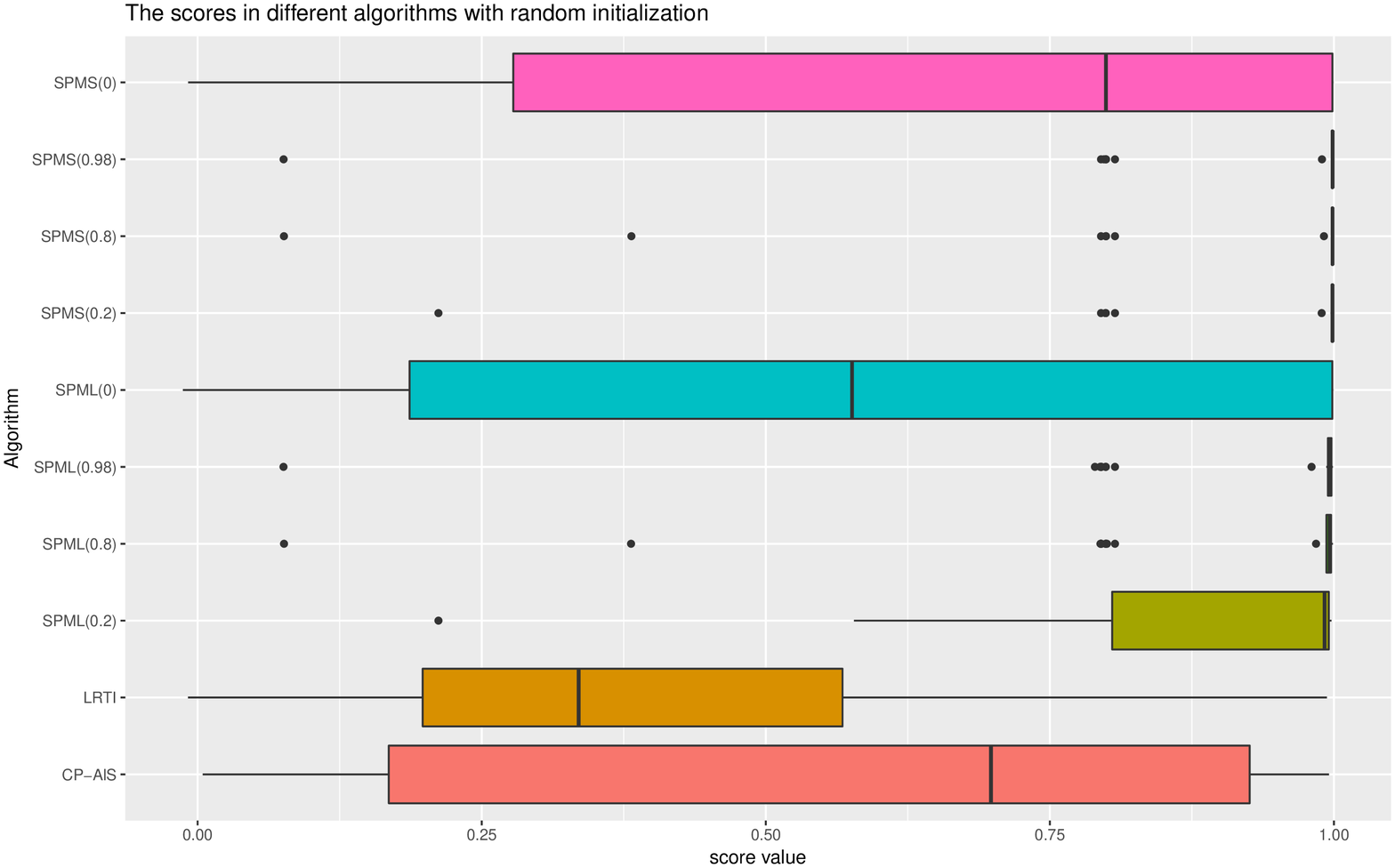}\\
\includegraphics[width=\columnwidth, height=0.6\textwidth ]{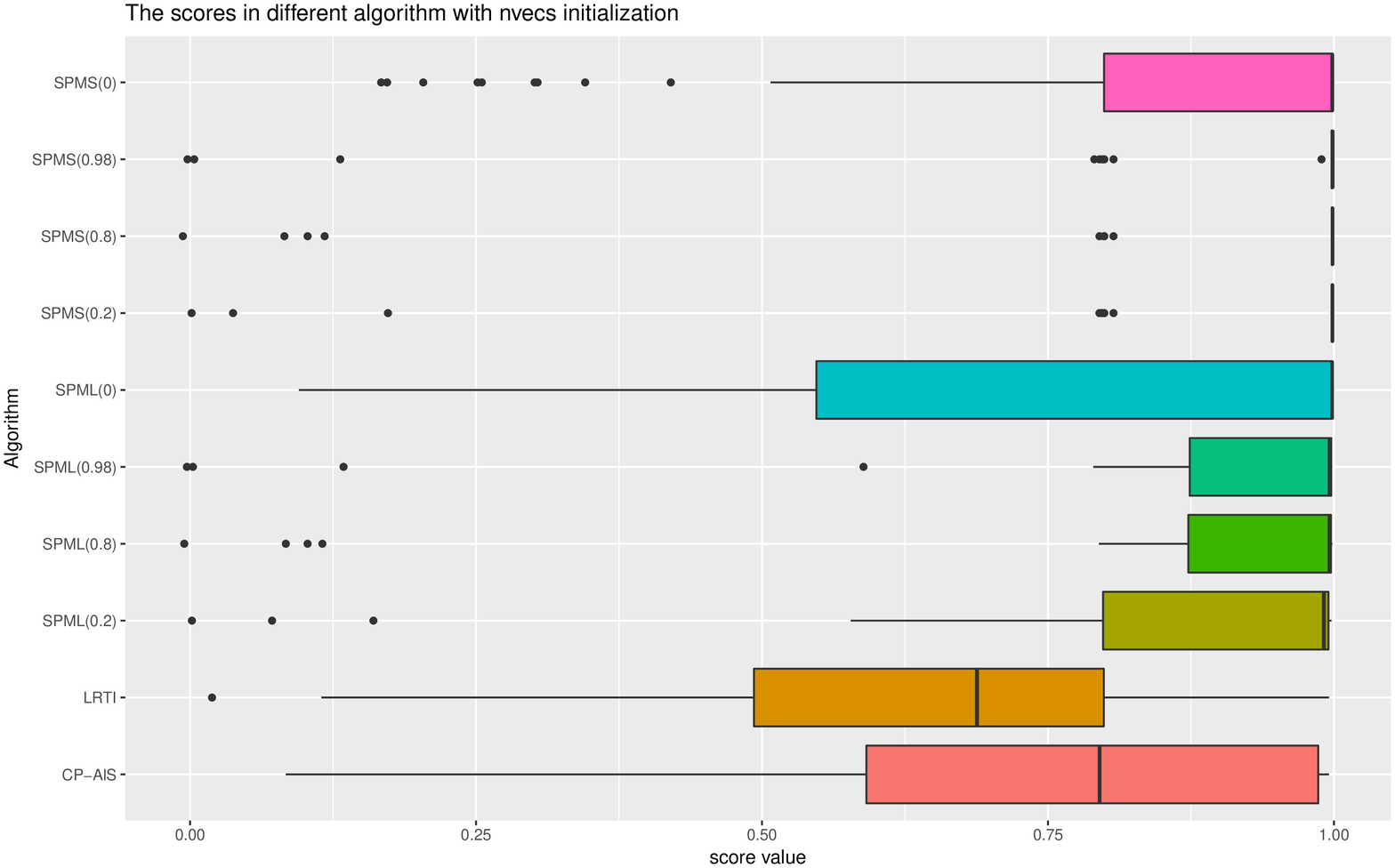}
\caption{The boxplot of scores for each algorithm in 50 random test when $(I_1, I_2, I_3)= (20,20,20)$, $R=10$  where the data with noise ( SNR = 20dB) and  without missing entries. The above figure using random initialization, and the below figure using nvecs initialzation.}
\label{fig:with-noise-no-missing20_10_10_random_scoresin_paper}
\end{center}
\vskip -0.2in
\end{figure}

And we also plot the scores for each algorithm with the random initialization in figure \ref{fig:with-noise-no-missing20_10_10_random_scoresin_paper}. 
%   \twocolumn[
%\newpage
%]

%]
%It can be seen that with an appropriate  guessing of $\alpha$, the proposed algorithm can find the true rank (M2). In most of the tests the true sparse factors can be approximated: when $\alpha = 0.8$ and in the sparse constrained model,  the ratio of total number of zeros in the selected solution and the total number of zeros in the true solution is 0.88 (M7), and the majority of zeros found by the algorithm are in the true solution (M9), sometimes the algorithm can even recover the true sparsity structure (M3, 6/50 times find the perfect sparse structure). When $\alpha = 0.8$,  the mean of total relative error in the selected solution subtracting the true relative error for Algorithm \ref{alg:lrti_coordinate_descent}  is 1.5e-2. The performance is further enhanced by the sparse constrained algorithm \ref{alg:lrti_coordinate_descent_sparse} to -8.3e-4. We see that the sparse constrained algorithm  \ref{alg:lrti_coordinate_descent_sparse}   is very useful to enhance the relative recovery error when the solution approximates the true sparsity structure. The convergence of algorithm\ref{alg:lrti_coordinate_descent}  is relatively slow (M15), and the sparse constrained algorithm \ref{alg:lrti_coordinate_descent_sparse}  converges much faster (M15).  
It can be seen that mean score value of our algorithm is better than the CP-ALS and LRTI algorithm. And we find our algorithm can find the sparse structure of the factor matrix. For  the mean of  the full data relative error   $\frac{||\underline{\bf Z} - \underline{\bf X} ||_F}{||\underline{\bf Z}||_F}$ subtract  the true relative error  $\frac{||\underline{\bf E}||_F}{||\underline{\bf Z}||_F}$, our method is as good as the CP-ALS method, and is better than the LRTI algorithm.  and SPMS(0.2) with random initialization achieves the best mean scores 0.967. And the most spare solution is find by the algorithm  SPMS(0.98) 
with random initialization. We find that the random initialization for our algorithm can achieve the better mean scores, while for the LRTI and CP-ALS algorithm, the nvecs initialization can achieve better mean scores. For the test time, that our algorithm may some costly, however, note that our algorithm calculate the whole solution path, the cost is relative small. 

\subsubsection{Model Performance with Noise and  Missing Values}
\label{subsubsection:random_test_behavior_noise_and_missing}
To gain more confidence for our mehthod, now we test on the synthetic data with noise and missing. Let $(I_1, I_2, I_3)= (50,50,50)$, $R=9$, the number of test times is set to $30$, and the initial estimate of rank is $10$.
We add noise to the data with the SNR of 20dB, where the SNR is defined by $\text{SNR} := 10\log_{10} \frac{\text{Var}(\underline{\bf X}) }{\text{Var}(\underline{\bf E})}$, and $\sigma^2 = \text{Var}(\underline{\bf E})$. And there are  missing values(random drop out 25\% entries). 
Since the CP-ALS algorithm cannot not hander the missing value entries, we only compare our method(both the solution path method(SPML, SPMS) and solution path method random(SPMLR , SPMSR)) with LRTI algorithm.

%\twocolumn[

The results are summarized in the table \ref{table:with-noise-with-missing50_9_10in_paper}. 

%\twocolumn[
%\newpage
%]

\begin{table}[!tb]
\caption{50 times random test behavior when $(I_1, I_2, I_3)= (50,50,50)$, $R=9$  where the data with noise ( SNR = 20dB) and  with missing entries(random drop out 25\% entries) and each algorithm using random initialization (given in above) and  nvecs initialization(given below).
And SPML stands for the Solution Path Method using  algorithm \ref{alg:lrti_coordinate_descent} (Look algorithm), and SPMS stands for the Solution Path Method using  algorithm \ref{alg:lrti_coordinate_descent_sparse} (Sparse constrained  Algorithm) where the solution is coming from  SPML. And SPML(0.8) standes for that $\alpha = 0.8$. And SCORES(STD) stand for the mean score and the its standard deviation. SPMLR stands for the Solution Path Method using  algorithm \ref{alg:lrti_coordinate_descent_subset},  And SPMSR stands for the Solution Path Method using  algorithm \ref{alg:lrti_coordinate_descent_sparse} (Sparse constrained  Algorithm) where the solution is coming from  SPMLR.
REL\_ERR stands for the mean of the full data relative error   $\frac{||\underline{\bf Z} - \underline{\bf X} ||_F}{||\underline{\bf Z}||_F}$ subtract  the true relative error  $\frac{||\underline{\bf E}||_F}{||\underline{\bf Z}||_F}$.
TFNZ stands for the total number of zeros found by the selected solution;  
TFTNZ stands for  total number of the zeros  such that the zero
appears in both the selected solution and the true solution;
TNZ stands for  total number of the zeros  in the true solution.
TIME stands for the the mean cost cpu time in 50 random testes, note that for each $\alpha$, SPML and SPMS are computed  together, so the time was spent in both the  SPML and SPMS in each test.
NFTR stands for number of find true rank in the 30 testes.
}
\label{table:with-noise-with-missing50_9_10in_paper}
\vskip 0.15in
\begin{center} 
%\begin{small}
\scriptsize
\begin{sc}
\begin{tabular}{cccccccccc}
\hline
\abovespace\belowspace
%$\lambda $&		N2 & N3 &N4 &N5&N6&N7\\
Method&		Scores(std) &rel\_err & TFNZ & TFTNZ & TNZ & $\frac{\text{TFNZ }}{\text{TNZ}}$ &$\frac{\text{TTNZ }}{\text{TNZ}}$   & time &NFTR\\
 \hline
\abovespace\belowspace

random &initialization \quad &\quad &\quad &\quad &\quad &\quad &\quad &\quad \\
\hline 
\abovespace
LRTI	&	0.726(0.256)	& 	6.641e-03	&  	     0	&  	     0	&  	  6985	&  	0.000	&  	0.000	& 	222.6&10\\ 
SPMLR(0.2)	&	0.851(0.154)	& 	1.008e-01	&  	   567	&  	   507	&  	  7064	&  	0.080	&  	0.072	& 	458.9 	&30\\ 
SPMSR(0.2)	&	1.000(0.000)	& 	-3.376e-04	&  	   567	&  	   567	&  	  7064	&  	0.080	&  	0.080	& 	458.9 	&30\\ 
SPMLR(0.8)	&	0.910(0.114)	& 	3.316e-02	&  	  2846	&  	  2619	&  	  7064	&  	0.403	&  	0.371	& 	494.6  	&30\\ 
SPMSR(0.8)	&	1.000(0.000)	& 	-3.200e-04	&  	  2846	&  	  2846	&  	  7064	&  	0.403	&  	0.403	& 	494.6  	&30\\ 
SPMLR(0.98)	&	0.957(0.093)	& 	1.216e-02	&  	  4658	&  	  4479	&  	  7064	&  	0.659	&  	0.634	& 	541.8  	&30\\ 
SPMSR(0.98)	&	1.000(0.000)	& 	-3.037e-04	&  	  4658	&  	  4658	&  	  7064	&  	0.659	&  	0.659	& 	541.8  	&30\\ 
SPML(0.2)	&	0.943(0.104)	& 	6.954e-02	&  	  6722	&  	  6454	&  	  7064	&  	0.952	&  	0.914	& 	322.1  	&30\\ 
SPMS(0.2)	&	1.000(0.000)	& 	-2.863e-04	&  	  6722	&  	  6722	&  	  7064	&  	0.952	&  	0.952	& 	322.1  	&30\\ 
SPML(0.8)	&	0.990(0.040)	& 	1.612e-02	&  	  7039	&  	  6992	&  	  7064	&  	0.996	&  	0.990	& 	322.2  	&30\\ 
SPMS(0.8)	&	1.000(0.000)	& 	-2.826e-04	&  	  7039	&  	  7039	&  	  7064	&  	0.996	&  	0.996	& 	322.2  	&30\\ 
SPML(0.98)	&	0.983(0.056)	& 	1.427e-02	&  	  7032	&  	  6948	&  	  7064	&  	0.995	&  	0.984	& 	323.1  	&30\\ 
\belowspace
SPMS(0.98)	&	1.000(0.000)	& 	-2.827e-04	&  	  7032	&  	  7032	&  	  7064	&  	0.995	&  	0.995	& 	323.1  	&30\\ 
\hline
\abovespace\belowspace
nvecs& initialization \quad \quad &\quad &\quad &\quad &\quad &\quad &\quad \\
\hline 
LRTI	&	0.670(0.236)	& 	7.102e-02	&  	     2	&  	     0	&  	  6985	&  	0.000	&  	0.000	& 	260.4&6\\ 
SPMLR(0.2)	&	0.836(0.158)	& 	1.067e-01	&  	   575	&  	   511	&  	  7064	&  	0.081	&  	0.072	& 	444.3  	&30\\ 
SPMSR(0.2)	&	1.000(0.000)	& 	-3.382e-04	&  	   575	&  	   575	&  	  7064	&  	0.081	&  	0.081	& 	444.3  	&30\\ 
SPMLR(0.8)	&	0.936(0.110)	& 	3.260e-02	&  	  2799	&  	  2649	&  	  7064	&  	0.396	&  	0.375	& 	487.3  	&30\\ 
SPMSR(0.8)	&	1.000(0.000)	& 	-3.205e-04	&  	  2799	&  	  2799	&  	  7064	&  	0.396	&  	0.396	& 	487.3  	&30\\ 
SPMLR(0.98)	&	0.954(0.106)	& 	1.160e-02	&  	  4307	&  	  4134	&  	  7064	&  	0.610	&  	0.585	& 	543.6  	&30\\ 
SPMSR(0.98)	&	1.000(0.000)	& 	-3.068e-04	&  	  4307	&  	  4307	&  	  7064	&  	0.610	&  	0.610	& 	543.6  	&30\\ 
SPML(0.2)	&	0.943(0.104)	& 	6.954e-02	&  	  6722	&  	  6454	&  	  7064	&  	0.952	&  	0.914	& 	657.3  	&30\\ 
SPMS(0.2)	&	1.000(0.000)	& 	-2.863e-04	&  	  6722	&  	  6722	&  	  7064	&  	0.952	&  	0.952	& 	657.3  	&30\\ 
SPML(0.8)	&	0.876(0.296)	& 	1.541e-02	&  	  6781	&  	  6196	&  	  7064	&  	0.960	&  	0.877	& 	610.1  	&26\\ 
SPMS(0.8)	&	0.886(0.297)	& 	-2.858e-04	&  	  6781	&  	  6243	&  	  7064	&  	0.960	&  	0.884	& 	610.1  	&26\\ 
SPML(0.98)	&	0.844(0.333)	& 	1.339e-02	&  	  7000	&  	  6138	&  	  7064	&  	0.991	&  	0.869	& 	590.4  	&26\\ 
\belowspace
SPMS(0.98)	&	0.868(0.336)	& 	-2.888e-04	&  	  7000	&  	  6274	&  	  7064	&  	0.991	&  	0.888	& 	590.4  	&26\\ 
\hline
\end{tabular}
\end{sc}
%\end{small}
\end{center}
\vskip -0.1in
\end{table}
%]
%   \twocolumn[
%\newpage
%]
% \twocolumn[

\begin{figure}[ht]
\vskip 0.2in
\begin{center}
\includegraphics[width=0.7\columnwidth,height=0.6\textwidth ]{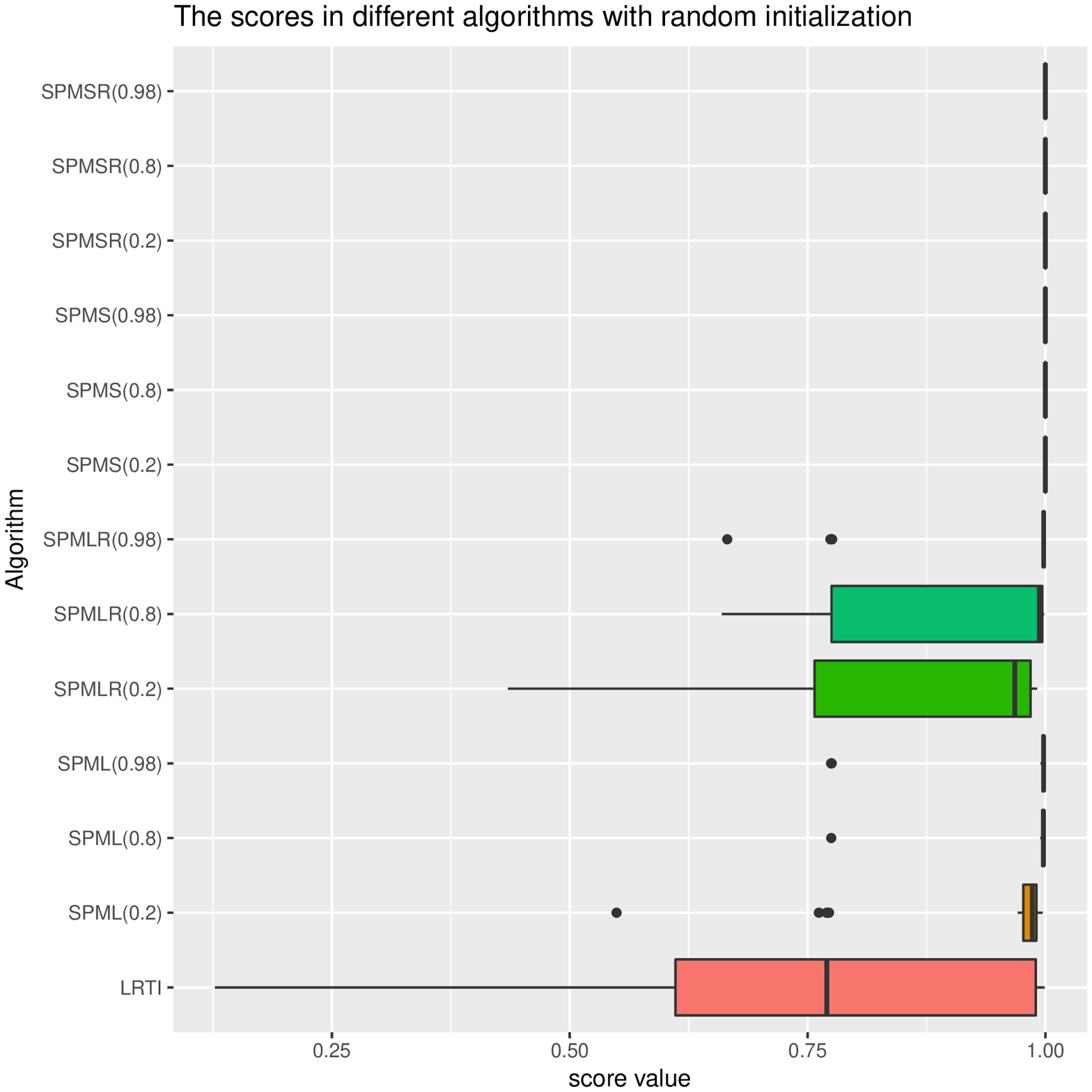}\\
\includegraphics[width=0.7\columnwidth,height =0.6\textwidth]{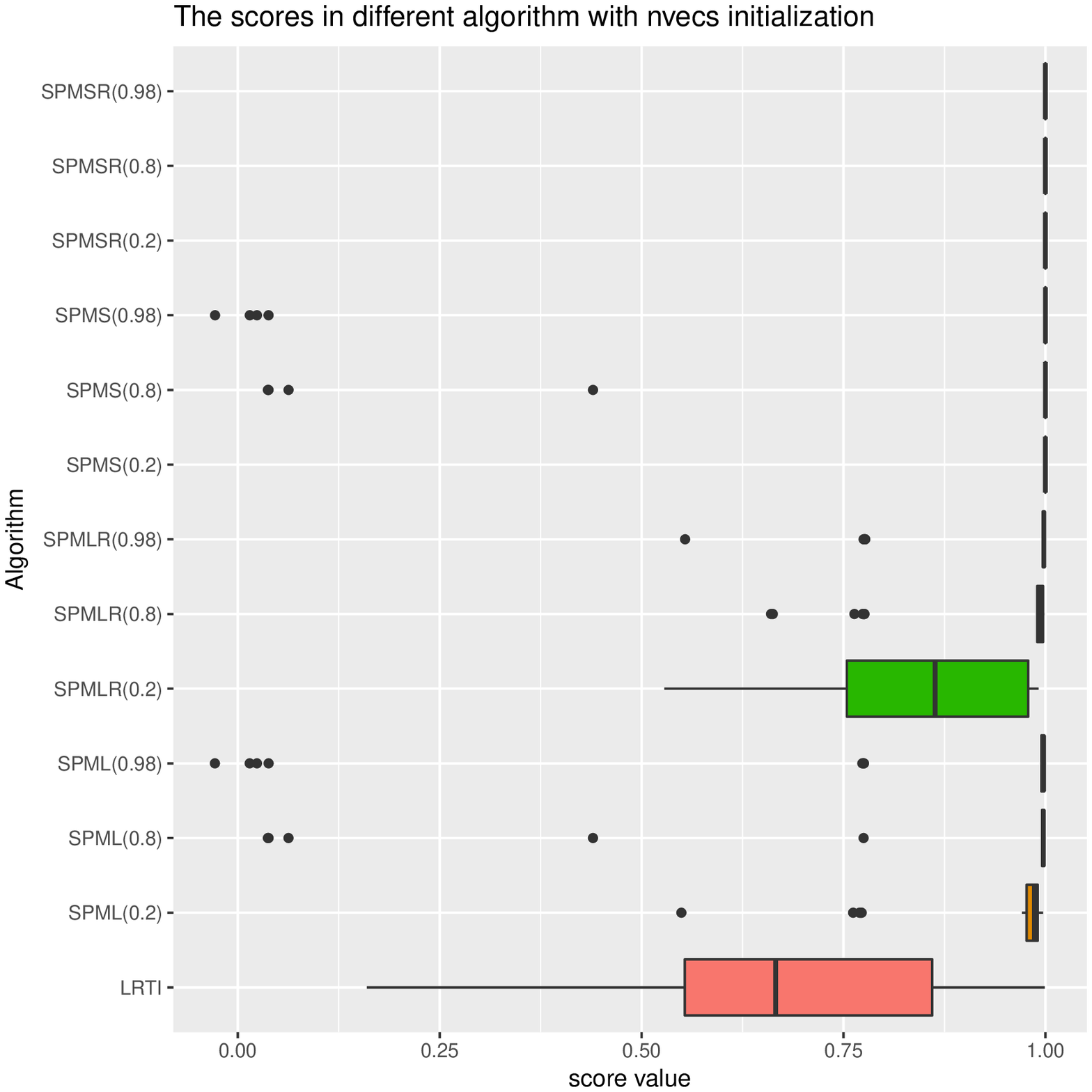}
\caption{The boxplot of scores for each algorithm  in 30 random test when $(I_1, I_2, I_3)= (50,50,50)$, $R=9$  where the data with noise ( SNR = 20dB) and  with random drop out 25\%  entries. The above figure using random initialization, and the below figure using nvecs initialzation.}
\label{fig:with-noise-with-missing50_9_10_scoresin_paper}
\end{center}
\vskip -0.2in
\end{figure}

And we also plot the scores for each algorithm in figure \ref{fig:with-noise-with-missing50_9_10_scoresin_paper}. 
%   \twocolumn[
%\newpage
%]

%]
%It can be seen that with an appropriate  guessing of $\alpha$, the proposed algorithm can find the true rank (M2). In most of the tests the true sparse factors can be approximated: when $\alpha = 0.8$ and in the sparse constrained model,  the ratio of total number of zeros in the selected solution and the total number of zeros in the true solution is 0.88 (M7), and the majority of zeros found by the algorithm are in the true solution (M9), sometimes the algorithm can even recover the true sparsity structure (M3, 6/50 times find the perfect sparse structure). When $\alpha = 0.8$,  the mean of total relative error in the selected solution subtracting the true relative error for Algorithm \ref{alg:lrti_coordinate_descent}  is 1.5e-2. The performance is further enhanced by the sparse constrained algorithm \ref{alg:lrti_coordinate_descent_sparse} to -8.3e-4. We see that the sparse constrained algorithm  \ref{alg:lrti_coordinate_descent_sparse}   is very useful to enhance the relative recovery error when the solution approximates the true sparsity structure. The convergence of algorithm\ref{alg:lrti_coordinate_descent}  is relatively slow (M15), and the sparse constrained algorithm \ref{alg:lrti_coordinate_descent_sparse}  converges much faster (M15).  
It can be seen that mean score value of our algorithm is better than the LRTI algorithm. And we find our algorithm can find the sparse structure of the factor matrix. For  the mean of  the full data relative error   $\frac{||\underline{\bf Z} - \underline{\bf X} ||_F}{||\underline{\bf Z}||_F}$ subtract  the true relative error  $\frac{||\underline{\bf E}||_F}{||\underline{\bf Z}||_F}$, our method is also better than the LRTI algorithm.  Many of our algorithm find out the true factors the score value are 1.0. And the random initialization is better than the nvecs initialization for both our algorithm and the LRTI algorithm. Most choices of our algorithm find out the true rank perfectly, while LRTI only find 10 times of true rank in 30 testes. For the test time, that our algorithm may some costly, however, note that our algorithm calculate the whole solution path, the cost is relative small.

\subsection{Application on Real Data}
\subsubsection{Coil-20 Data}
\begin{figure}[ht]
\vskip 0.2in
\begin{center}
\includegraphics[width=0.7\columnwidth,height=0.6\textwidth ]{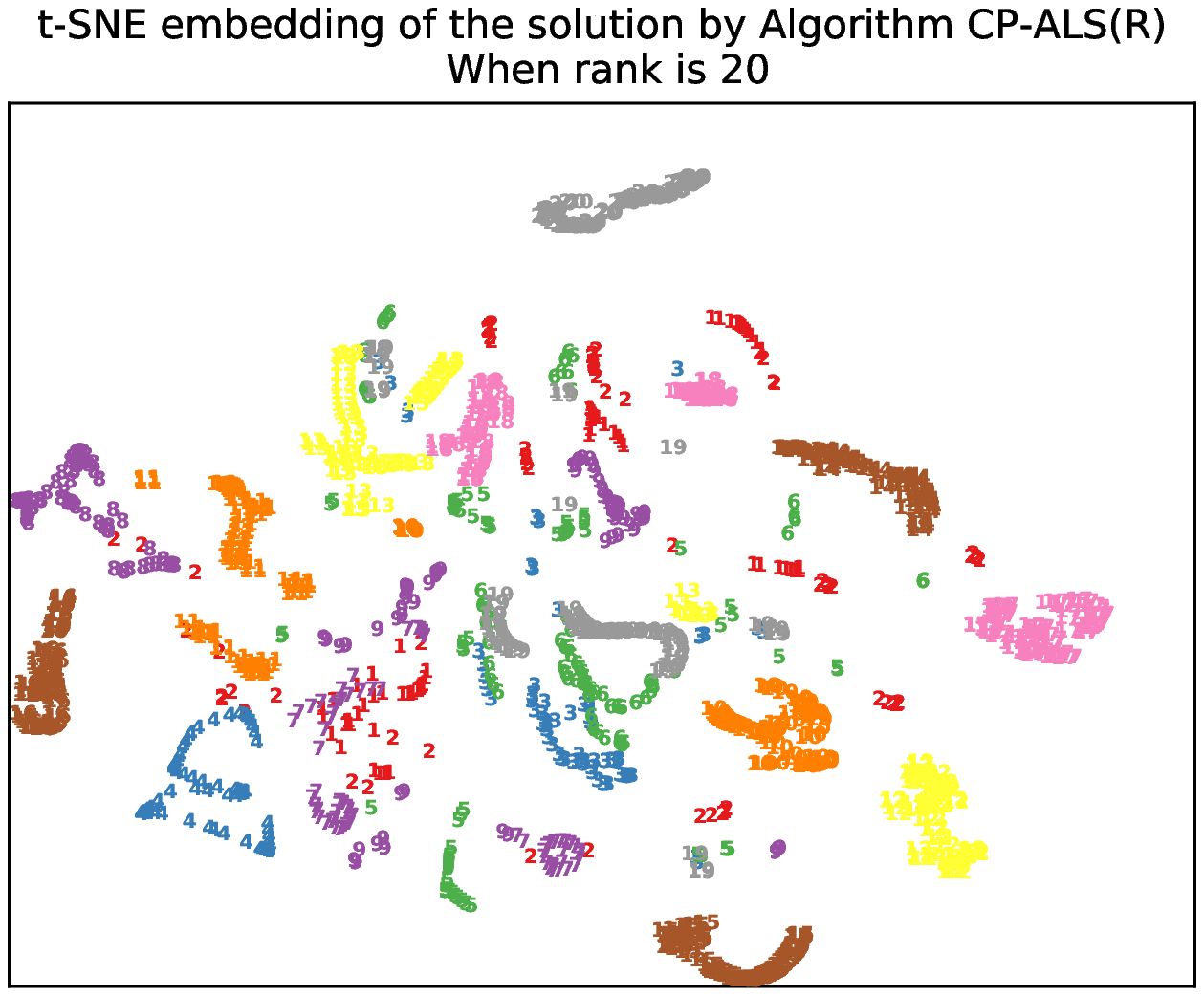}\\
\includegraphics[width=0.7\columnwidth,height =0.6\textwidth]{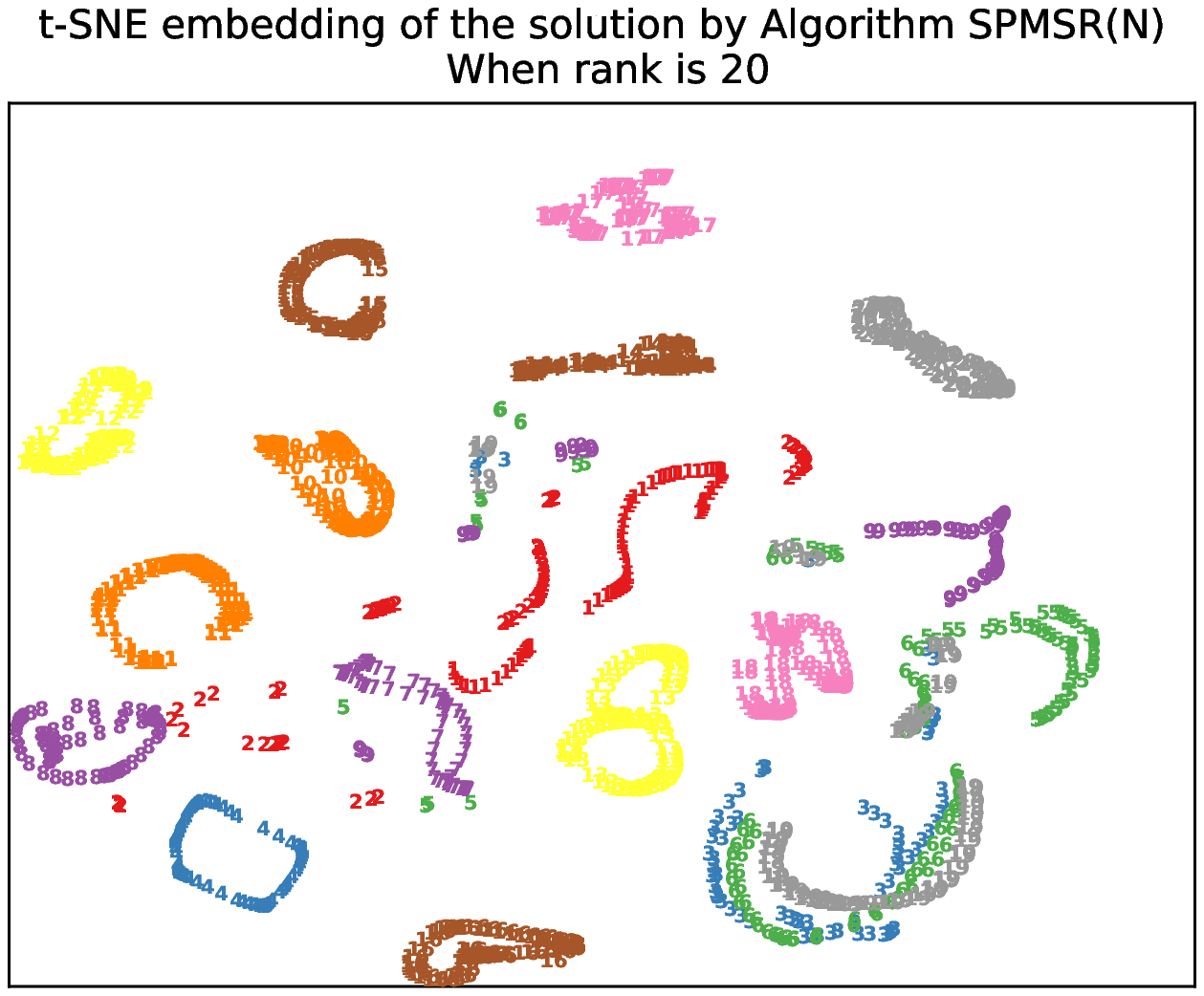}
\caption{Performance of  CP-ALS(R) and SPMSR(0.8N) when initial rank is 20 in the COIL-20 image clustering. The number in the figure stand the object index. The two t-SNE components of ${\bf A}^{(3)}$ were used for visualization.}
\label{fig:coil-20_rank_20}
\end{center}
\vskip -0.2in
\end{figure}

\begin{figure}[ht]
\vskip 0.2in
\begin{center}
\includegraphics[width=0.7\columnwidth,height=0.6\textwidth ]{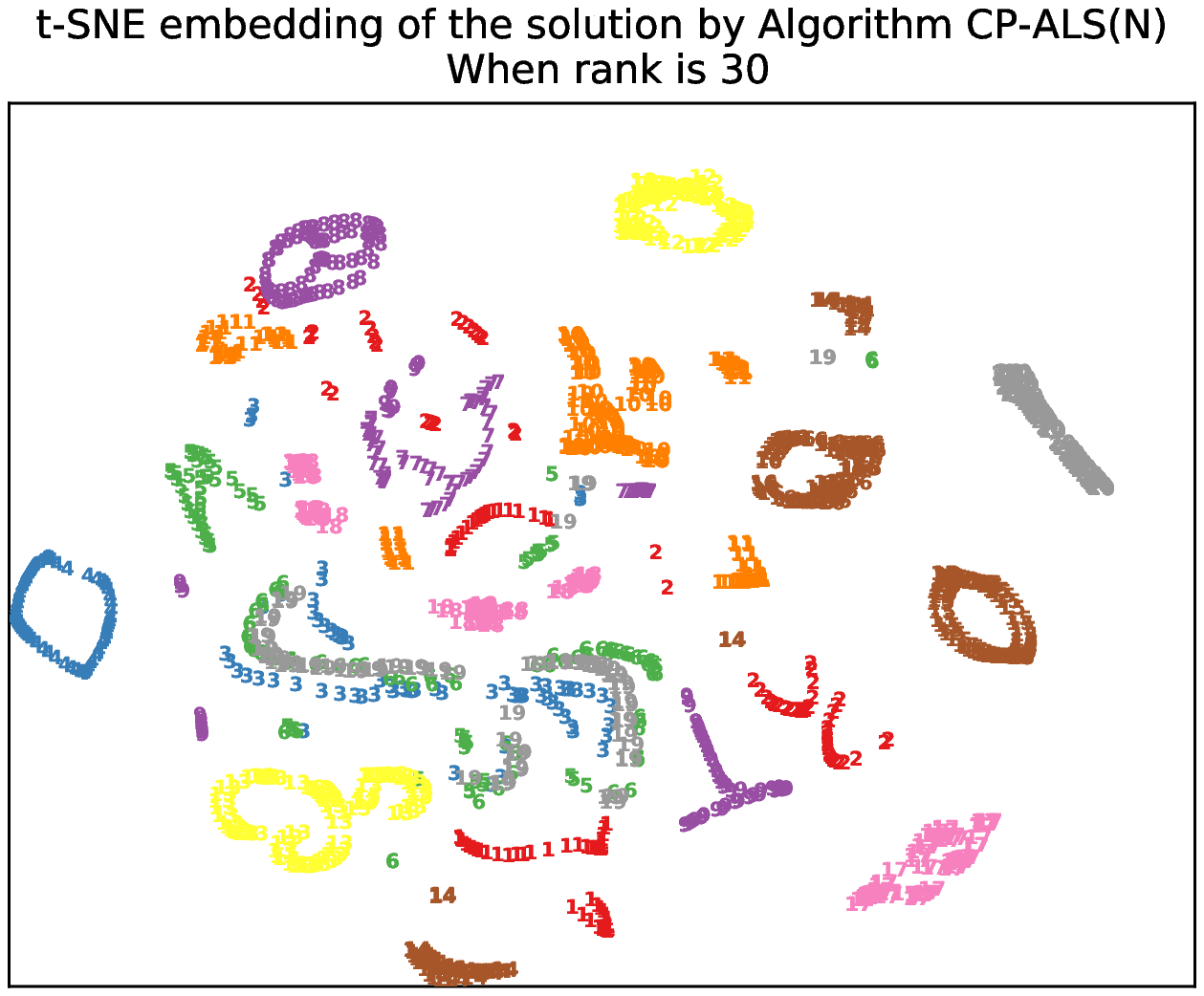}\\
\includegraphics[width=0.7\columnwidth,height =0.6\textwidth]{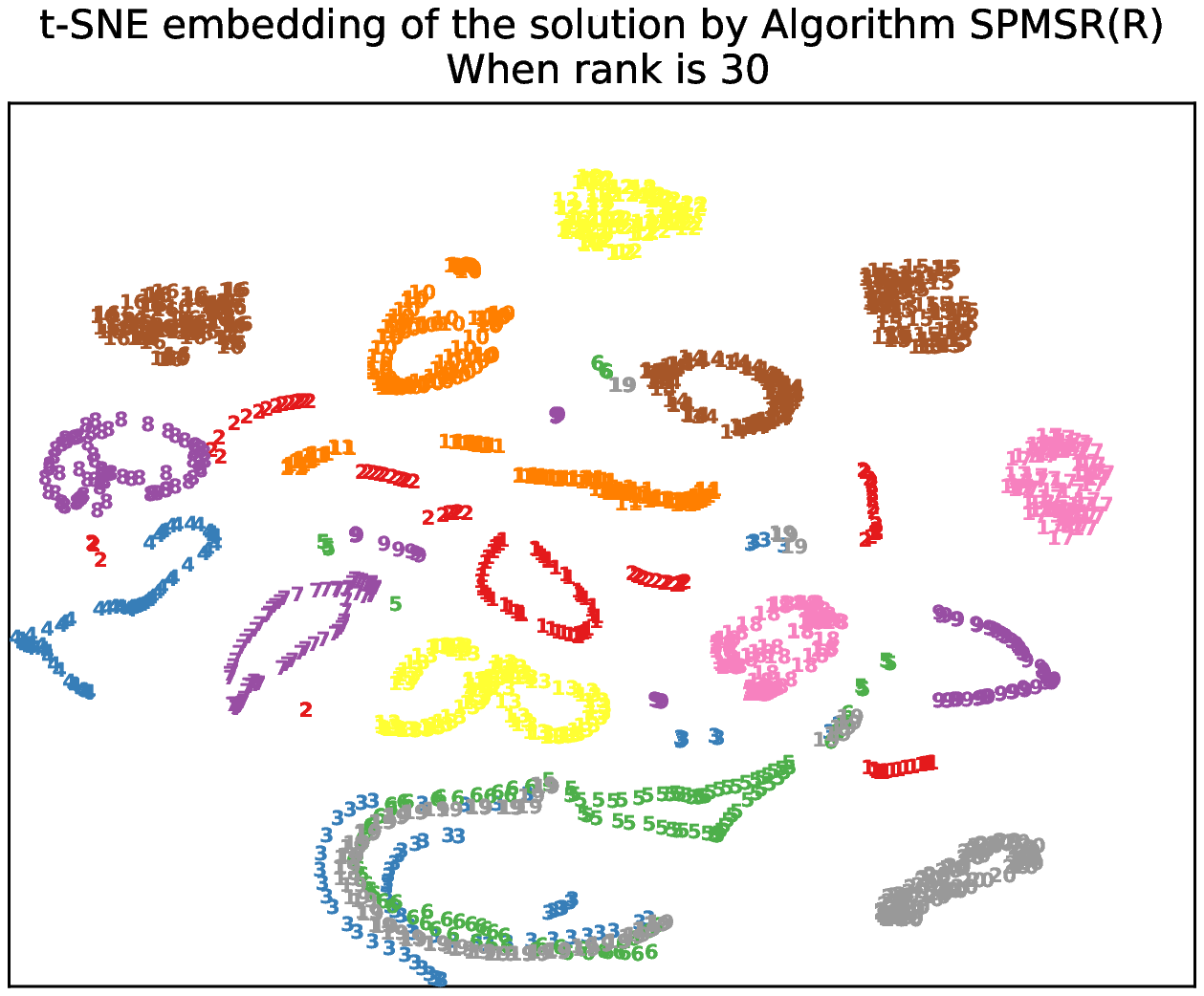}
\caption{Performance of  CP-ALS(N) and SPMSR(0.8R) when initial rank is 30 in the COIL-20 image clustering.The number in the figure stand the object index. The two t-SNE components of ${\bf A}^{(3)}$ were used for visualization.}
\label{fig:coil-20_rank_30}
\end{center}
\vskip -0.2in
\end{figure}

In this section, we test the proposed algorithm on the coil-20 data (\cite{coil-20}).  This database contains
1440 images of 20 objects, each image is $128 \times 128$ pixels, and each object is captured from
varying angles with a 5-degree interval. We generated
a tensor $\underline{\bf Z}$ with $128 \times 128 \times 1440$ using these images.  And we compare  our algorithms  (SPMLR, SPMSR) with the CP-ALS and LRTI algorithm(Not given the result of the LRTI algorithm, since it behaves poor, the relative error is more than $0.9$). We don't use the   (SPML, SPMS) because it a little slow and the performance of SPMLR, SPMSR are as good as SPML, SPMS. 
And we use the two t-SNE (\cite{t-SNE})  components of ${\bf A}^{(3)}$ for visualization and
clustering. The K-means algorithm was adopted for clustering. As
K-means is prone to be affected by initial cluster centers, in each run we repeated clustering 20 times, each with the init 'k-means++'(see for the help for K-means method in scikit-learn(\cite{scikit-learn}) ). The performance averaged over 20 Monte-Carlo
runs is detailed in Table \ref{table:coil-20}. We give the t-SNE visualize of CP-ALS(R) and SPMSR(0.8N) when initial rank is 20 in Figure \ref{fig:coil-20_rank_20}, and the t-SNE visualize of CP-ALS(N) and SPMSR(0.8R) when initial rank is 20 in Figure \ref{fig:coil-20_rank_30}. From these results, we find that algorithm is good for find the underline true factors, although the relative recovery error is a little high than the the CP-ALS method, however, it can achieve a better cluster result (from both the table and the figure) than the CP-ALS method. Which that algorithm is trend to find the true factors.

\begin{table}[h]
\caption{Performance comparison of the algorithms in COIL-20 image
clustering, where Rank is the initial rank, where R and N in the first row in the parenthesis are stand for the random initialization and nvecs initialization respectively.}
\label{table:coil-20}
\vskip 0.15in
\begin{center} 
%\begin{small}
\scriptsize
\begin{sc}
\begin{tabular}{cccccccccc}
\hline
\abovespace\belowspace
%$\lambda $&		N2 & N3 &N4 &N5&N6&N7\\
Algorithms &CP-ALS(R) &  CP-ALS(N)& SPMLR(0.8R)&SPMSR(0.8R) )& SPMLR(0.8N)&SPMSR(0.8N) \\
 \hline
 \abovespace\belowspace
 Rank=20\\ \hline
 \abovespace\belowspace
 
RelErr & 		0.294&	0.291&0.316	&0.301	&   	0.309 &	0.299&&\\
\abovespace
\belowspace
NumZeros&		0&		0&	31 &31	&		26&		26& \\
\abovespace
\belowspace
accuracy(\%)& 	0.637&   0.602&	0.700&	0.724&	0.760&	0.762&&&\\
\abovespace
\belowspace
time(s) &		514.8&	478.1&40882.5&	40882.5&	45547.4&45547.4&\\
 \hline
 \abovespace\belowspace
 Rank=30 \\\hline
 \abovespace\belowspace
 RelErr & 		0.265&	0.263&	0.296&	0.278&	0.284&	0.272&\\
 \abovespace
\belowspace
NumZeros&		0&		0&		30&	30&		23&		23& \\
\abovespace
\belowspace
\abovespace
\belowspace
accuracy(\%) &	0.561&	0.670&	0.717 &	0.752 &	0.676&	0.751&&\\
\abovespace
\belowspace
time(s) &		550.7&	566.8&	43123.7&	43123.7&	47380.1&47380.1&\\

\hline
\end{tabular}
\end{sc}
%\end{small}
\end{center}
\vskip -0.1in
\end{table}

\subsubsection{COPD Data}
In this work, we test the proposed algorithm on the baseline data from The Genetic Epidemiology of COPD (COPDGene$^{\textregistered}$) cohort (\cite{copd}). The data record consists of 10,300 subjects and 361 features for each subject, including self-administered questionnaires of demographic data and medical history, symptoms, medical record review, etc. All features are independently reviewed by certified professionals (the full data collection forms can be found www.COPDGene.org). 

For the data preprocessing, we firstly remove the features which cannot be quantified, then discard the features with missing data covering more than 20 percent of the whole record, thus finally form a table of $213$ features. We then normalize the data by $\frac{{\bf x} - \min({\bf x})}{\max({\bf x}) - \min({\bf x})}$ for each feature ${\bf x}$ except the "GOLD value" feature which grades the COPD stages. As there are $7$ GOLD values( $-2, -1, \ldots , 4$) for $7$ stages and the number of subjects within each stage is different (minimum is $108$), we perform subsampling on each stage for balancing the data.  

To reveal the features which are closely related to the severity of COPD, we transform the 2D (subject by feature) table into a 3-way tensor $\underline{\bf Z}$, whose dimension is $(7, 108, 212)$, where the first dimension is the GOLD values, the second dimension is the subjects, and the third dimension is the features. We also keep records for the location of the missing values in tensor $\underline{ \boldsymbol \Delta}$. 

\begin{figure}[ht]
\vskip 0.2in
\begin{center}
\centerline{\includegraphics[width=\columnwidth]{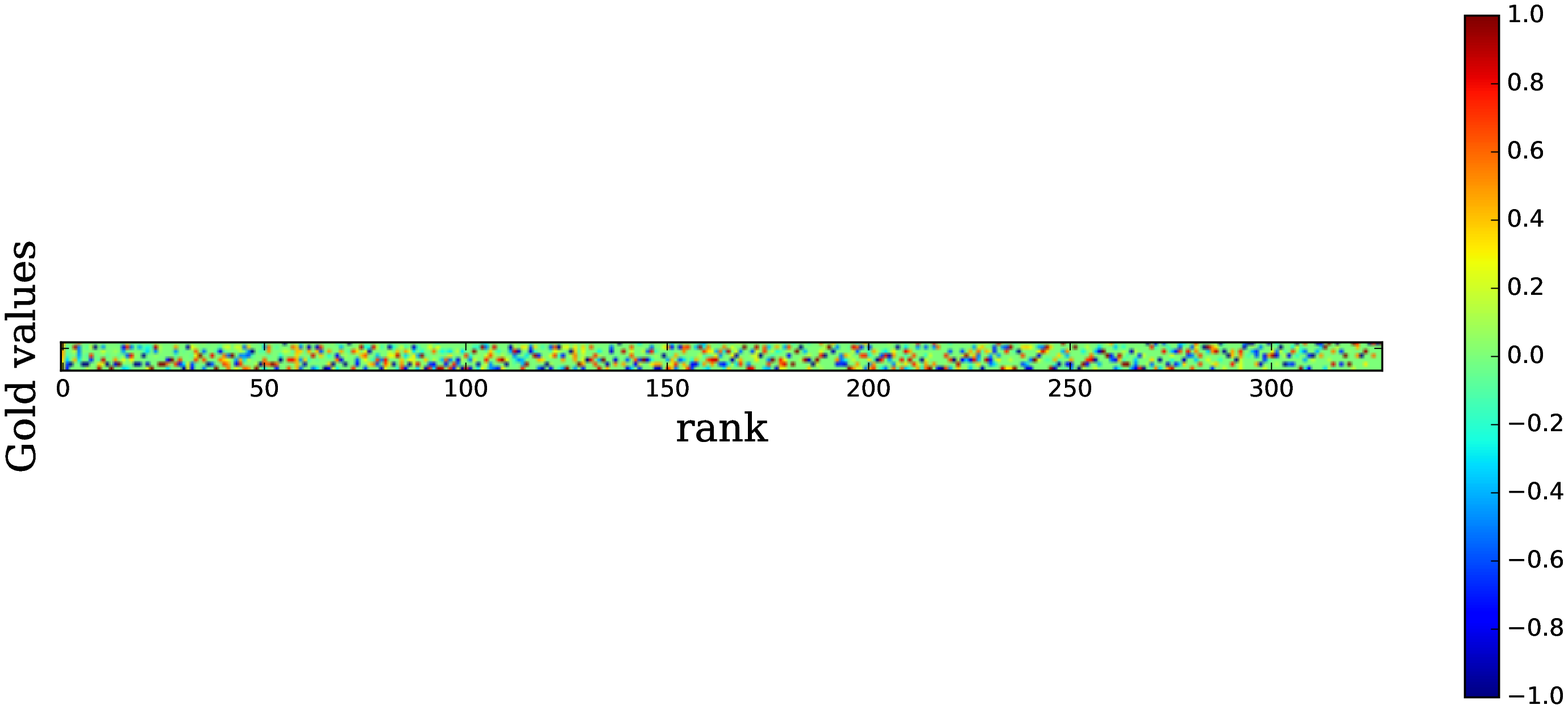}}
\caption{The figure of the first normalized  factor matrix and for the Gold values}
\label{fig:factor_matrix0}
\end{center}
\vskip -0.2in
\end{figure}

\begin{figure}[ht]
\vskip 0.2in
\begin{center}
\centerline{\includegraphics[width=\columnwidth]{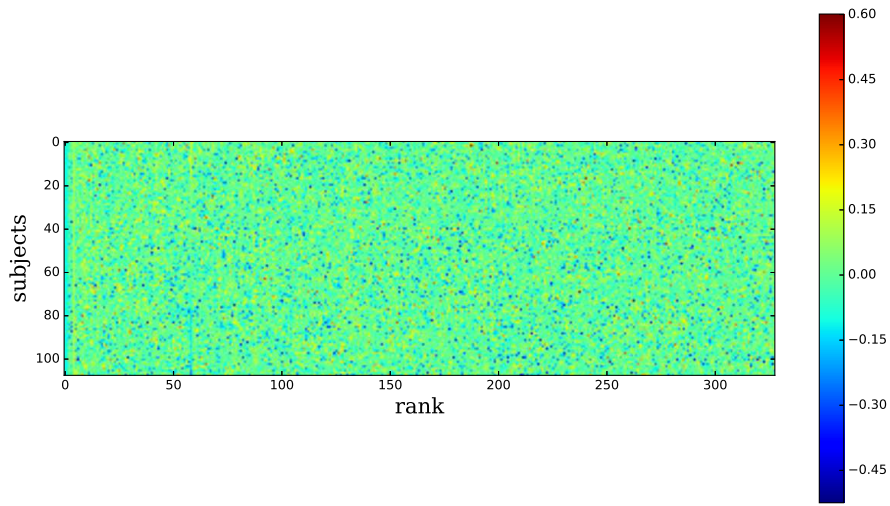}}
\caption{The figure of the second normalized  factor matrix for the subjects.}
\label{fig:factor_matrix1}
\end{center}
\vskip -0.2in
\end{figure}

\begin{figure}[ht]
\vskip 0.2in
\begin{center}
\centerline{\includegraphics[width=\columnwidth]{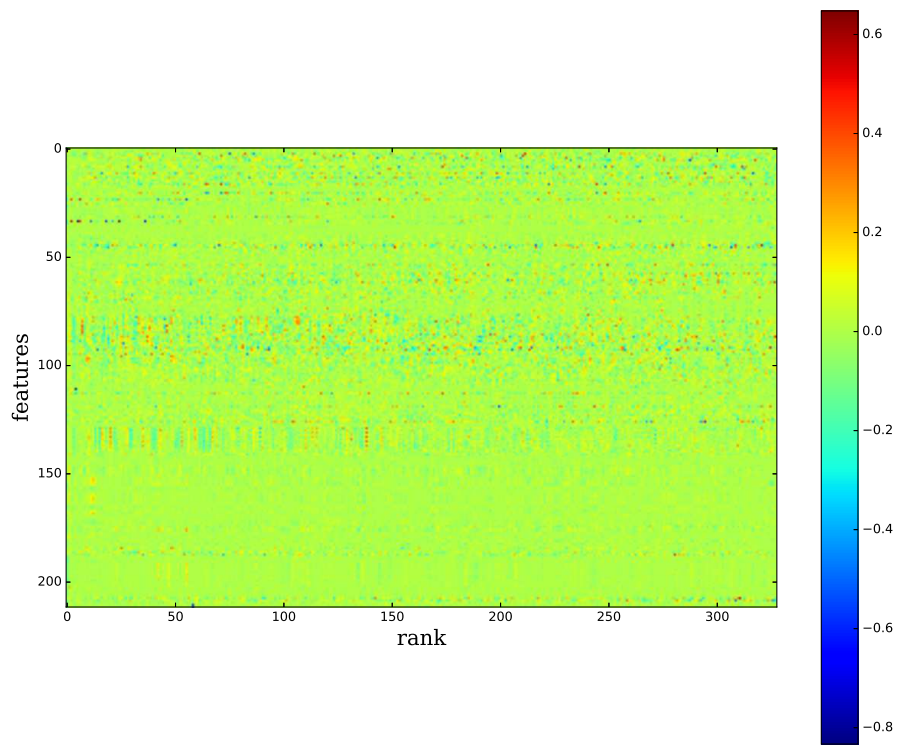}}
\caption{The figure of the third normalized factor matrix for the features.}
\label{fig:factor_matrix2}
\end{center}
\vskip -0.2in
\end{figure}

\begin{figure}[ht]
\vskip 0.2in
\begin{center}
\centerline{\includegraphics[width=\columnwidth]{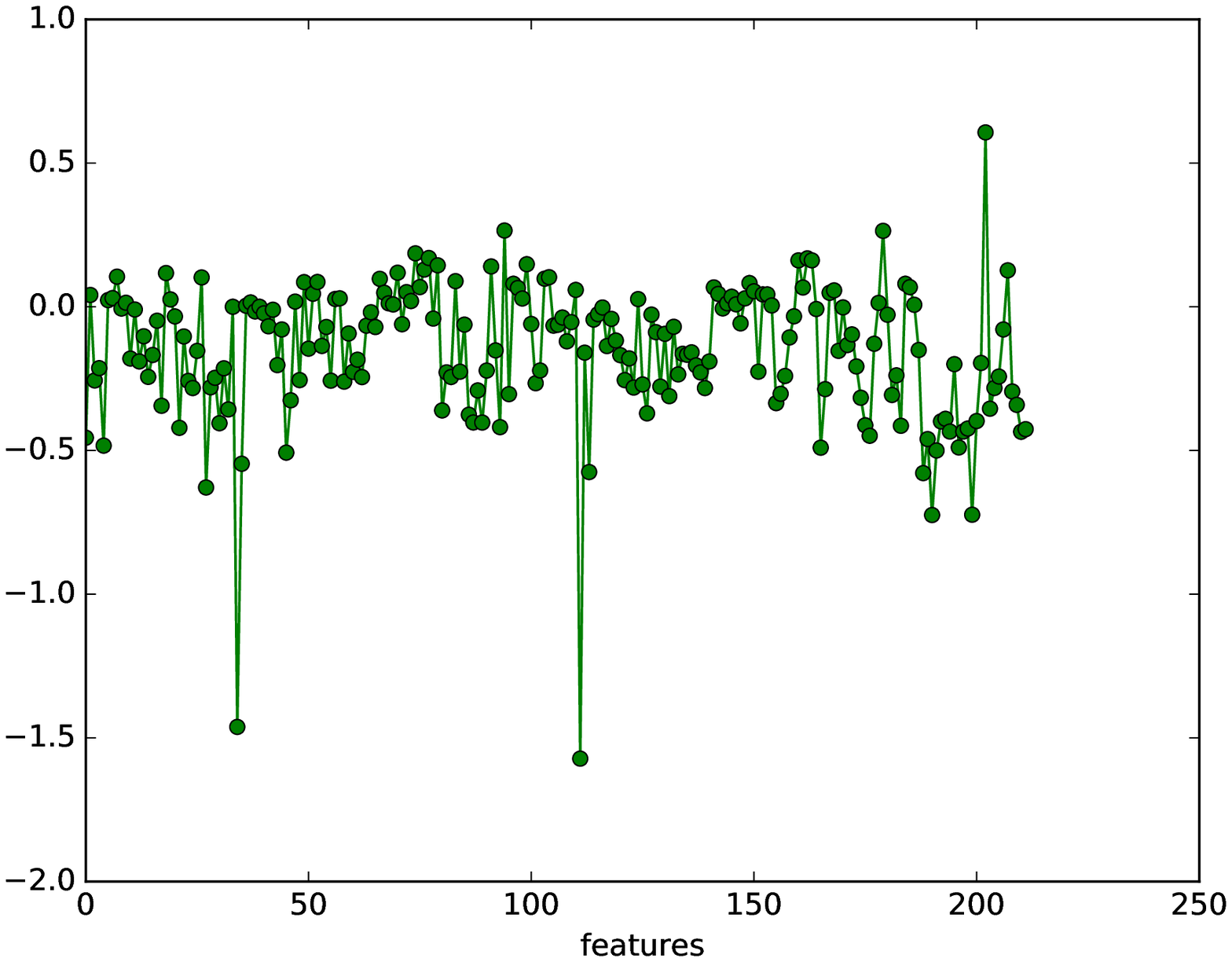}}
\caption{The figure of the sum of the first $6$ columns (corresponding to the largest $6$ character values) of the third feature factor matrix. The largest top 10 features(sorted by the absolute value) are 111(EverSmokedCig): -1.6,
    34(Blood\_Other\_Use): -1.5,
   190(FEV1\_FVC\_utah): -0.7,    199(pre\_FEV1\_FVC): -0.7,
    27(Resting\_SaO2): -0.6,
   202(ATS\_ERS): 0.6,    188(FEV1pp\_utah): -0.6,
   113(SmokCigNow): -0.6,
    35(HealthStatus): -0.5,     45(HighBloodPres): -0.5.
%And we print the top 10 absolute value in the title
}
\label{fig:row_sum}
\end{center}
\vskip -0.2in
\end{figure}

By applying Algorithm \ref{alg:solution_path} on the 3-way tensor above, we generate a solution path and select the best solution from it. The rank of the best solution is $328$ and the relative recovery error is $0.137$. The sparseness (the number of zeros elements found / total number of elements in the factor matrices) is $0.26$. As shown in the visualization of the three factor matrices in Figure
\ref{fig:factor_matrix0}, \ref{fig:factor_matrix1}, \ref{fig:factor_matrix2}, all of them have sparse structure. We further calculate the sum of the top $6$ columns (which correspond the largest $6$ character values in ${\boldsymbol \gamma}$) of the third features's factor matrix   to get a vector of $212$ elements and sort them by their absolute value (see in Figure \ref{fig:row_sum}). The top $10$ features are EverSmokedCig, Blood\_Other\_Use, FEV1\_FVC\_utah, pre\_FEV1\_FVC, Resting\_SaO2, ATS\_ERS, FEV1pp\_utah, SmokCigNow, HealthStatus,  HighBloodPres. According to domain experts' experience, most of these features are closely related the severity of COPD.

%
%The figure of the row sum of the frontal $6$ columns (which correspond to the largest 6 character values in ${\boldsymbol \gamma}$)  of the third feature factor matrix is given in Figure \ref{fig:row_sum}. 

\section{Conclusion}
\label{sec:conclusion}
In this work, we derive a tensor decomposition model from the Bayesian framework to learn the sparse and low rank PARAFAC decomposition via elastic net regularization. An efficient block coordinate descent algorithm and stochastic block coordinate descent algorithm are proposed 
to solve the optimization problem,
%for this problem,
 which can be applied to solve large scale problems. To increase the robustness of the algorithm for practical applications, we develop a solution path strategy with warm start to generate a gradually sparse and low rank solution, and we reduce the relative recovery error by using the sparse constrained algorithm \ref{alg:lrti_coordinate_descent_sparse}. Evaluation on  synthetic data shows our algorithm can capture the true rank and 
 find a close approximation of the true
% closed to
  true sparse structure with high probability and   with the help of sparse structure, the constrained algorithm trends to  find the true underline factors. Evaluation on the coil-20 data shows that our algorithm can extract the meaningful underline factors that can have better cluster performance. Evaluation on the real   data shows that it can extract the sparse factors that reveal meaningful relationships between the features in the data and the severity of COPD.

In the future, we plan to expand the current algorithm to the N-way nonnegative PARAFAC decomposition in a similar structure.
%
%\section{Conclusion}
%\label{sec:conclusion}
%In this article, we derive a model from the Bayesian framework to learning the sparse and low rank PARAFAC decomposition via the elastic net regularization for the observed tensor with missing values. An efficient coordinate descent algorithm is proposed for this problem, which can be used to solve large scale problem. To choose the reasonable solution in practical, we use the solution path with warm start strategy to generate a gradually sparse and low rank solution path, and we enhance the relative recovery error by using the sparse constrained algorithm \ref{alg:lrti_coordinate_descent_sparse}. From the numerical test on the synthetic data, our algorithm can capture the true rank and closed to true sparse structure with high probability. And when applied on the COPD data, our algorithm extracts the sparse factor matrix which contains meaningful information about the severity of COPD.
%
%We will plan to enhance the speed and save the storage of our algorithm by devising a stochastic coordinate descent scheme. Also we plan to  design a similar version for the N-way nonnegative PARAFAC decomposition in the future work.
%
% In the unusual situation where you want a paper to appear in the
% references without citing it in the main text, use \nocite
\nocite{reg_path}
\nocite{fast_cd}
\nocite{sntf_ccd}
\nocite{learning_parts}
%\nocite{fast_cd}

\newpage
\bibliography{refer.bib}

\begin{thebibliography}{16}
\providecommand{\natexlab}[1]{#1}
\providecommand{\url}[1]{\texttt{#1}}
\expandafter\ifx\csname urlstyle\endcsname\relax
  \providecommand{\doi}[1]{doi: #1}\else
  \providecommand{\doi}{doi: \begingroup \urlstyle{rm}\Url}\fi

\bibitem[Bazerque et~al.(2013)Bazerque, Mateos, and Giannakis]{rankreg}
Juan~Andres Bazerque, Gonzalo Mateos, and Georgios~B. Giannakis.
\newblock Rank regularization and bayesian inference for tensor completion and
  extrapolation.
\newblock \emph{IEEE Transactions on Signal Processing}, 61\penalty0
  (22):\penalty0 5689--5703, 2013.

\bibitem[Chartrand(2007)]{lp-sparse-solution}
R.~Chartrand.
\newblock Exact reconstruction of sparse signals via nonconvex minimization.
\newblock \emph{IEEE Signal Processing Letters}, 14\penalty0 (10):\penalty0
  707--710, Oct 2007.
\newblock ISSN 1070-9908.
\newblock \doi{10.1109/LSP.2007.898300}.

\bibitem[Cichocki and Phan(2009)]{fast_cd}
Andrzej Cichocki and Anh~Huy Phan.
\newblock Fast local algorithms for large scale nonnegative matrix and tensor
  factorizations.
\newblock \emph{Ieice Trans Fundamentals}, 92\penalty0 (3):\penalty0 708--721,
  2009.

\bibitem[Friedman et~al.(2010)Friedman, Hastie, and Tibshirani]{reg_path}
Jerome Friedman, Trevor Hastie, and Rob Tibshirani.
\newblock Regularization paths for generalized linear models via coordinate
  descent.
\newblock \emph{Journal of Statistical Software}, 33\penalty0 (i01):\penalty0
  1, 2010.

\bibitem[Hinton(2008)]{t-SNE}
G.~E. Hinton.
\newblock Visualizing high-dimensional data using t-sne.
\newblock \emph{Vigiliae Christianae}, 9\penalty0 (2):\penalty0 2579--2605,
  2008.

\bibitem[Kingma and Ba(2014)]{adam}
Diederik~P Kingma and Jimmy Ba.
\newblock Adam: A method for stochastic optimization.
\newblock \emph{Computer Science}, 2014.

\bibitem[Kolda and Bader(2006)]{tensortoolbox}
Tamara~G Kolda and Brett~W Bader.
\newblock Matlab tensor toolbox.
\newblock 2006.

\bibitem[Kolda and Bader(2009)]{td_review}
Tamara~G Kolda and Brett~W Bader.
\newblock Tensor decompositions and applications.
\newblock \emph{Siam Review}, 51\penalty0 (3):\penalty0 455--500, 2009.

\bibitem[Lee and Seung(1999)]{learning_parts}
D.~D. Lee and H.~S. Seung.
\newblock Learning the parts of objects by non-negative matrix factorization.
\newblock \emph{Nature}, 401\penalty0 (6755):\penalty0 788--791, 1999.

\bibitem[Liu et~al.(2012)Liu, Liu, Wonka, and Ye]{sntf_ccd}
Ji~Liu, Jun Liu, Peter Wonka, and Jieping Ye.
\newblock Sparse non-negative tensor factorization using columnwise coordinate
  descent.
\newblock \emph{Pattern Recognition}, 45\penalty0 (1):\penalty0 649--656, 2012.

\bibitem[Pedregosa et~al.(2011)Pedregosa, Varoquaux, Gramfort, Michel, Thirion,
  Grisel, Blondel, Prettenhofer, Weiss, Dubourg, Vanderplas, Passos,
  Cournapeau, Brucher, Perrot, and Duchesnay]{scikit-learn}
F.~Pedregosa, G.~Varoquaux, A.~Gramfort, V.~Michel, B.~Thirion, O.~Grisel,
  M.~Blondel, P.~Prettenhofer, R.~Weiss, V.~Dubourg, J.~Vanderplas, A.~Passos,
  D.~Cournapeau, M.~Brucher, M.~Perrot, and E.~Duchesnay.
\newblock Scikit-learn: Machine learning in {P}ython.
\newblock \emph{Journal of Machine Learning Research}, 12:\penalty0 2825--2830,
  2011.

\bibitem[Regan et~al.(2010)Regan, Hokanson, Murphy, Make, Lynch, Beaty,
  Curraneverett, Silverman, and Crapo]{copd}
E.~A. Regan, J.~E. Hokanson, J.~R. Murphy, B~Make, D.~A. Lynch, T.~H. Beaty,
  D~Curraneverett, E.~K. Silverman, and J.~D. Crapo.
\newblock Genetic epidemiology of copd (copdgene) study design.
\newblock \emph{Copd-journal of Chronic Obstructive Pulmonary Disease},
  7\penalty0 (1):\penalty0 32, 2010.

\bibitem[S.~A.~Nene and H.(1996)]{coil-20}
S.~K.~Nayar S.~A.~Nene and Murase H.
\newblock Columbia object image library (coil-20).
\newblock \emph{Technical Report CUCS-005-96}, February 1996.

\bibitem[Tomasi and Bro(2006)]{comparealg}
Giorgio Tomasi and Rasmus Bro.
\newblock A comparison of algorithms for fitting the parafac model.
\newblock \emph{Computational Statistics and Data Analysis}, 50\penalty0
  (7):\penalty0 1700--1734, 2006.

\bibitem[Xu and Yin(2013)]{bcd}
Yangyang Xu and Wotao Yin.
\newblock A block coordinate descent method for regularized multiconvex
  optimization with applications to nonnegative tensor factorization and
  completion.
\newblock \emph{Siam Journal on Imaging Sciences}, 6\penalty0 (3):\penalty0
  1758--1789, 2013.

\bibitem[Zou and Hastie(2005)]{elastic_net}
Hui Zou and Trevor Hastie.
\newblock Regularization and variable selection via the elastic net.
\newblock \emph{Journal of the Royal Statistical Society}, 67\penalty0
  (2):\penalty0 301--320, 2005.

\end{thebibliography}

\newpage

\section{Appendix}
\subsection{Proof of Theorem \ref{th:convergence}}
\label{supp:theorem}
%\begin{theorem}
%%[\ref{th:convergence}]
%Let $\{ {\bf A}^k\} $ be the sequence generated by Algorithm  \ref{alg:lrti_coordinate_descent}, where  $\{ {\bf A}^k\} $  is the solution $({\bf A}^{(1)}, {\bf A}^{(2)}, \ldots , {\bf A}^{(N)} )$ in  the k-th iteration in the repeat loop. Assume that $\{ {\bf A}^k\} $ is bound. Then $\{ {\bf A}^k\} $ converges to  a critical point $\bar {\bf A}$, and the asymptotic convergence rates in Theory 2.9 in \cite{bcd} apply.
%\end{theorem*}
%
\begin{proof}
From equation (\ref{obj:bcd:f}),   (\ref{obj:bcd:r}),  (\ref{obj:bcd:F}), our notations coincide with the notation in  (\cite{bcd}), (where  in our algorithm, the block $i$ is just a column ${\bf a}^{(n)}_{r}$, for $r = 1, 2, \ldots, R, n=1, 2, \ldots, N)$. Since our algorithm \ref{alg:lrti_coordinate_descent} is a special case of Algorithm 1 in  (\cite{bcd}) in which we only use the update rule (1.3a) in (\cite{bcd})
, we only need to verify the assumptions in Theory 2.8 and Theory 2.9 in (\cite{bcd}).

For the assumption 1  in  (\cite{bcd}): Obviously, $F$ is continuous in $\text{dom}(F)$ and $\inf_{\bf{A}} F(\bf{A}) \ge 0$, and $ F(\bf{A})$ has a Nash point (See (2.3) in  (\cite{bcd}) for definition).

For the assumption 2  in  (\cite{bcd}) :   
\begin{equation}\label{obj:bcd:f_element}
\begin{aligned}
f({\bf a}^{(n)}_{r}) := &\frac{1}{2 } || (\underline{\bf Z} - \sum_{r=1}^R {\bf a}_r^{(1)} \circ {\bf a}_r^{(2)} \circ \cdots \circ {\bf a}_r^{(N)}) {\varoast} \underline{{\boldsymbol {\Delta}}}  ||_F^2 \\
				&  + \sum_{r= 1}^R \sum_{n= 1}^N   \lambda\left [     \frac{1-\alpha}{2} (({\bf a}^{(n)}_r)^T {\bf R}_n^{-1}{\bf a}^{(n)}_r)  \right ]\\
 \end{aligned}
\end{equation}       
 From the equation (\ref{subgradient:F}), it is easy to see that $\nabla^2_{{\bf a}^{(n)}_{r}} f( {\bf a}^{(n)}_{r}  ) = {\bf H} + \lambda (1-\alpha) {\bf R}_n^{-1}$, where {\bf H} is defined by (\ref{H}) with ${\bf h} = {\bf a}_r^{{\varoast}_{-n}} :={\bf a}^{(N)}_r \varotimes \cdots \varotimes  {\bf a}^{(n+1)}_r \varotimes  {\bf a}^{(n-1)}_r \varotimes \cdots \varotimes  {\bf a}^{(1)}_r$, and ${\boldsymbol{\delta}}_{i_n}^T$ is the $i_n$-th row of ${\boldsymbol{\Delta}}_{(n)}$. So $f({\bf a}^{(n)}_{r}) $ is strongly convex with modulus $l {\bf I}\le \nabla^2_{{\bf A}^{(n)}_{i_n,r}} f(  {\bf A}^{(n)}_{i_n,r}  )  \le L {\bf I}$, where $l:= \lambda(1-\alpha) \underset{n=1,2, \ldots, N, i_n = 1, 2, \ldots, I_N}{\min} {\bf R}_n^{-1}(i_n, i_n) >0 $, since $\alpha = \frac{\mu}{1+\mu} <1$. And $L$  exist since we assume that $\{ {\bf A}^k\} $ is bounded. 
 
 For the conditions in Lemma 2.6  in  (\cite{bcd}):
 \begin{enumerate}
\item   From the verification of the assumption 2  in  (\cite{bcd}) above, it is easy to see that $\nabla f$ is Lipschitz continuous on any bounded set
\item $F$ satisfies the KL inequility (2.14) in  (\cite{bcd}) at $\bar {\bf A}$ since $f$  satisfies the KL inequility (2.14) and $r$  satisfies KL inequility (2.14), so does their sum.
\item We may choose the initial estimate ${\bf A}$ sufficiently close to $\bar{\bf A}$, and $F({\bf A}^k) >F(\bar{\bf A})$ for $k>0$, since $F({\bf a}^{(n)}_{r}) $ is strongly convex in ${\bf a}^{(n)}_{r}$, so $F({\bf A}^k)$ is strictly decreasing. 
\end{enumerate}
All the conditions in  Theory 2.8 and Theory 2.9 in (\cite{bcd}) are satisfied, so the conclusions follow out.
\end{proof}

\subsection{Proof of Proposition \ref{pro:l2equi}}
\label{supp:prop}
\begin{proof}
Let $\tilde {\bf A} := (\tilde {\bf A}^{(1)}, \tilde {\bf A}^{(2)}, \ldots,  \tilde {\bf A}^{(N)})$, and $\tilde {\bf A}^{(n)} = {\bf R}_n^{-1/2} {\bf A}^{(n)}$, $\tilde {\bf a}_r^{(n)} := {\bf R}_n^{-1/2} {\bf a}_r^{(n)}$, then we can rewrite (\ref{tensor_imputation_with_l2_norm}) as 
 \begin{equation}\label{tensor_imputation_with_l2_norm_standard_sub}
 \begin{aligned}
\underset{\tilde {\bf A}, \;\underline{\bf X} }{\min} &\;\frac{1}{2 } || (\underline{\bf Z} - \underline{\bf X}) {\varoast} \underline{{\boldsymbol {\Delta}}}  ||_F^2   + \sum_{r= 1}^R \sum_{n= 1}^N       \frac{\lambda}{2} ||\tilde{\bf a}^{(n)}_r||_2^2 \\
			\text{s.t.}&\;	 \underline{\bf X} = \sum_{r=1}^R  ({\bf R}_1^{1/2} \tilde {\bf a}_r^{(1)} )   \circ \cdots \circ ( {\bf R}_N^{1/2} \tilde {\bf a}_r^{(N)})
 \end{aligned}
\end{equation}
We can standardize $\tilde {\bf a}_r^{(n)}$ to unit length, $\tilde {\bf u}_r^{(n)}:= \frac{\tilde {\bf a}_r^{(n)}}{||\tilde {\bf a}_r^{(n)}||_2}$, $\tilde {\bf U}^{(n)}:= (\tilde {\bf u}_r^{(1)}, \tilde {\bf u}_r^{(2)}, \ldots  , \tilde {\bf u}_r^{(N)})$,  $\tilde {\gamma}_r := \prod_{n=1}^N ||\tilde {\bf a}_r^{(n)}||_2, \; r = 1, 2, \ldots R$, and let ${a}_r^{(n)}:= ||\tilde {\bf a}_r^{(n)}||_2$. 
Then (\ref{tensor_imputation_with_l2_norm_standard_sub}) is equivalent to
 \begin{equation}\label{tensor_imputation_with_l2_norm_standard_2_sub}
 \begin{aligned}
	&\underset{\tilde {\bf U}^{(n)}, \; \tilde \gamma_r , \; \tilde{a}^{(n)}_r,\; \underline{\bf X}}{ \min} \frac{1}{2 } || \left (\underline{\bf Z} - \underline{\bf X} \right ) {\varoast} \underline{{\boldsymbol {\Delta}}}  ||_F^2 	+ \sum_{r= 1}^R  \frac{\lambda}{2} \sum_{n= 1}^N       (\tilde{a}^{(n)}_r)^2  \\
		& \text{s.t. }\tilde \gamma_r = \prod_{n=1}^N \tilde{a}^{(n)}_r \; , r = 1, 2, \ldots , R\\
	& \underline{\bf X} =\sum_{r=1}^R  \tilde \gamma_r ({\bf R}_1^{1/2} \tilde {\bf u}_r^{(1)} )  \circ \cdots \circ ( {\bf R}_N^{1/2} \tilde {\bf u}_r^{(N)})
\end{aligned}
\end{equation}
% ,\; n= 1, 2,\ldots ,N ,\; r = 1, 2, \ldots , R
Focus on the inner minimization $\text{w.r.t.}$ norms $ \{ \tilde{a}^{(n)}_r \}$ for arbitrary directions $ \{ \tilde{\bf a}^{(n)}_r \}$ and fixed products $\tilde {\gamma}_r$ and $\underline{\bf X}$ , the (\ref{tensor_imputation_with_l2_norm_standard_2_sub}) is equivalent to
 \begin{equation}\label{tensor_imputation_with_l2_norm_standard_2_inner_supp}
 \begin{aligned}
\underset{ {\tilde{a}^{(n)}_r, \; n= 1, 2,\ldots ,N  }} {\min} & \sum_{n= 1}^N       (\tilde{a}^{(n)}_r)^2  \\
		\text{s.t. } & \tilde \gamma_r = \prod_{n=1}^N \tilde{a}^{(n)}_r
\end{aligned}
\end{equation}
The arithmetic geometric-mean inequality gives the solution to (\ref{tensor_imputation_with_l2_norm_standard_2_inner_supp}), as it states that for scalars $\{ \tilde{a}^{(n)}_r \}$, it holds that
\begin{equation}\label{geometric-mean inequality}
 \begin{aligned}
(\prod_{n=1}^N(\tilde{a}^{(n)}_r)^2  )^{\frac{1}{N}} \le \frac{1}{N} \sum_{n=1}^N (\tilde{a}^{(n)}_r)^2  
\end{aligned}
\end{equation}
with equality when $\tilde{a}^{(1)}_r  = \tilde{a}^{(2)}_r = \cdots =\tilde{a}^{(N)}_r $ so that the minimum of (\ref{tensor_imputation_with_l2_norm_standard_2_inner_supp}) is attained at  $\tilde{a}^{(1)}_r  = \tilde{a}^{(2)}_r = \cdots =\tilde{a}^{(N)}_r   = ( {\tilde {\gamma}_r} )^{\frac{1}{N}}$.

Then (\ref{tensor_imputation_with_l2_norm_standard_2_sub}) is equivalent to
 \begin{equation}\label{tensor_imputation_with_l2_norm_standard_3_supp}
 \begin{aligned}
&\underset{\tilde {\bf U}^{(n)}, \; \tilde \gamma_r ,\; \underline{\bf X}}{ \min} \frac{1}{2 } || (\underline{\bf Z} - \underline{\bf X}) {\varoast} \underline{{\boldsymbol {\Delta}}}  ||_F^2  + \sum_{r= 1}^R  \frac{\lambda N }{2}       ( {\tilde {\gamma}_r} )^{\frac{2}{N} } \\
				& \text{s.t. } \underline{\bf X} =\sum_{r=1}^R  \tilde \gamma_r ({\bf R}_1^{1/2} \tilde {\bf u}_r^{(1)} )  \circ \cdots \circ ( {\bf R}_N^{1/2} \tilde {\bf u}_r^{(N)})
\end{aligned}
\end{equation}
So the solution of (\ref{tensor_imputation_with_l2_norm}) and (\ref{tensor_imputation_with_l2_norm_standard_3_supp}) coincide, i.e. the optimal factors related by  
$ {\bf a}_r^{(N)} = (\tilde {\gamma}_r)^{\frac{1}{N}} {\bf R}_n^{1/2}\tilde {\bf u}_r^{(n)} $
\end{proof}

\subsection{The relation of kernel similarity matrix and covariance matrix}
Note that ${\bf X}_{(n)} = {\bf A}_n( {\bf A}_N \varodot \cdots \varodot  {\bf A}_{n+1} \varodot  {\bf A}_{n-1} \varodot \cdots \varodot  {\bf A}_{1})^T$, let ${\bf B}_n :=  {\bf A}_N \varodot \cdots \varodot  {\bf A}_{n+1} \varodot  {\bf A}_{n-1} \varodot \cdots \varodot  {\bf A}_{1}$, then
\begin{equation}\label{kernel_similarity_matrix_1}
 \begin{aligned}
&{\bf X}_{(n)} (i,:) ({\bf X}_{(n)} (j,:))^T  = \; {\bf A}_n(i,:){\bf B}_n( {\bf B}_n )^T({\bf A}_n(j,:))^T \\
			 = &\; {\bf A}_n(i,:) [  (({\bf A}_{1})^T{\bf A}_{1}) \varoast  \cdots \varoast (({\bf A}_{n-1})^T  {\bf A}_{n-1}) \varoast \\
			& \; (({\bf A}_{n+1})^T{\bf A}_{n+1}) \varoast  \cdots   \varoast (({\bf A}_{N})^T{\bf A}_{N}) ]  ({\bf A}_n(j,:))^T  \\
			=&\;   \text{Tr}( [  (({\bf A}_{1})^T{\bf A}_{1}) \varoast  \cdots \varoast (({\bf A}_{n-1})^T  {\bf A}_{n-1})  \varoast  \\
			& \; (({\bf A}_{n+1})^T{\bf A}_{n+1}) \varoast  \cdots   \varoast (({\bf A}_{N})^T{\bf A}_{N}) ]\\
			&[({\bf A}_n(j,:))^T  {\bf A}_n(i,:)]  )
 \end{aligned}
\end{equation} 
Let $\theta_n := \text{Tr}({\bf R}_n)$.
Then 
\begin{equation}\label{R:estimate1}
 \begin{aligned}
 \mathbb E (({\bf A}_{n})^T{\bf A}_{n}) = \theta_n {\bf I}_R
 \end{aligned}
\end{equation} 
\begin{equation}\label{R:estimate2}
 \begin{aligned}
 \mathbb E [({\bf A}_n(j,:))^T  {\bf A}_n(i,:)]  ={\bf R}_n(i,j)  {\bf I}_R
 \end{aligned}
\end{equation} 
Combine equation (\ref{kernel_similarity_matrix} ), (\ref{kernel_similarity_matrix_1} ), (\ref{R:estimate1} ), (\ref{R:estimate2} ), we get 
\begin{equation}\label{kernel_similarity_matrix_sol}
{\bf K}_n = R \prod_{m=1, m \ne n}^N \theta_m {\bf R}_n
\end{equation} 
and 
\begin{equation}\label{X_F_norm}
\mathbb E ||\underline{\bf X}||_F^2 = \text{Tr}({\bf K}_n)= R \prod_{m=1}^N \theta_m
\end{equation}
Since from our assumption (\ref{equi_power}), $\theta_n = \theta, \; n=1, 2, \ldots, N$, we get
\begin{equation}\label{kernel_similarity_matrix_sol2}
{\bf K}_n = R\theta^{N-1} {\bf R}_n
\end{equation} 
and 
\begin{equation}\label{X_F_norm2}
\mathbb E ||\underline{\bf X}||_F^2 = \text{Tr}({\bf K}_n)=  R\theta^{N}
\end{equation}

\subsection{Full Solution Path}
\label{supp:solution_path}
We give the full solution path in  table \ref{table:solution path full}.

%
%
%\begin{table}[!tb]
%\caption{Solution path by Algorithm \ref{alg:solution_path} on the simultated  data when $\alpha = 0.2$. R: rank, NZS: the number of zeros in the current solution, NZT: the number of the  zeros such that the zero appears in the current solution and the true solution;
%IS1: which algorithm is used?  1--algorithm \ref{alg:lrti_coordinate_descent} and 0-- algorithm \ref{alg:lrti_coordinate_descent_sparse}; 
%}
%\vskip 0.1in
%\label{table:solution path}
%\begin{center} 
%\begin{small}
%\begin{sc}
%\begin{tabular}{ccccccc}
%\hline
%\abovespace\belowspace
%%$\lambda $&		N2 & N3 &N4 &N5&N6&N7\\
%$\lambda $&		R & NZS &NZT &IS1&rel\_err&iters\\
% \hline
% \abovespace

\begin{table}[ht]
\caption{Solution path by Algorithm \ref{alg:solution_path} on the simultated  data when $\alpha = 0.2$. R: rank, NZS: the number of zeros in the current solution, NZT: the number of the  zeros such that the zero appears in the current solution and the true solution;
IS1: which algorithm is used?  1--algorithm \ref{alg:lrti_coordinate_descent} and 0-- algorithm \ref{alg:lrti_coordinate_descent_sparse}; 
}
\vskip 0.15in
\label{table:solution path full}
\begin{center} 
%\begin{small}
\scriptsize
\begin{sc}
\begin{tabular}{ccccccc}
\hline
\abovespace\belowspace
$\lambda $&		R & NZS &NZT &IS1&rel\_err&iters\\
 \hline
 \abovespace
1e-10&     3&     3&     2& 1&8.8e-06& 38 \\  
1e-08&     3&     3&     2& 0&6.4e-06& 2 \\  
1e-09&     3&     3&     2& 1&6.4e-06& 2 \\  
1e-08&     3&     3&     2& 1&5.3e-06& 1 \\  
1e-07&     3&     3&     2& 1&4.4e-06& 1 \\  
1e-06&     3&     3&     2& 1&2.4e-06& 3 \\  
1e-05&     3&     4&     2& 1&2.1e-06& 3 \\  
1e-08&     3&     4&     2& 0&1.6e-06& 1 \\  
1e-04&     3&     3&     2& 1&1.8e-05& 56 \\  
1e-08&     3&     3&     2& 0&3.5e-06& 3 \\  
1e-03&     3&     3&     2& 1&1.7e-04& 200 \\  
1e-02&     3&     3&     2& 1&1.2e-03& 200 \\  
1e-01&     3&     8&     4& 1&3.7e-03& 200 \\  
1e-08&     3&     8&     4& 0&2.3e-05& 100 \\  
1e+00&     2&     5&     5& 1&2.3e-03& 200 \\  
1e-08&     2&     5&     5& 0&4.4e-06& 12 \\  
1e+01&     2&     7&     7& 1&1.3e-02& 200 \\  
1e-08&     2&     7&     7& 0&3.2e-06& 15 \\  
2e+01&     2&     8&     8& 1&2.5e-02& 200 \\  
1e-08&     2&     8&     8& 0&4.0e-06& 16 \\  
3e+01&     2&     8&     8& 1&3.6e-02& 200 \\  
4e+01&     2&     8&     8& 1&4.8e-02& 200 \\  
5e+01&     2&     8&     8& 1&5.8e-02& 200 \\  
6e+01&     2&     8&     8& 1&6.9e-02& 200 \\  
7e+01&     2&     8&     8& 1&7.9e-02& 178 \\  
8e+01&     2&     8&     8& 1&8.8e-02& 160 \\  
9e+01&     2&     8&     8& 1&9.8e-02& 146 \\  
1e+02&     2&     8&     8& 1&1.1e-01& 134 \\  
2e+02&     2&    10&    10& 1&1.9e-01& 108 \\  
1e-08&     2&    10&    10& 0&3.3e-06& 17 \\  
3e+02&     2&    10&    10& 1&2.6e-01& 63 \\  
4e+02&     2&    11&    11& 1&3.2e-01& 63 \\  
4e+02&     2&    11&    11& 0&4.4e-06& 8 \\  
5e+02&     2&    12&    12& 1&3.8e-01& 53 \\  
1e-08&     2&    12&    12& 0&1.6e-06& 9 \\  
6e+02&     2&    12&    12& 1&4.4e-01& 45 \\  
7e+02&     1&     3&     3& 1&6.1e-01& 56 \\  
1e-08&     1&     3&     3& 0&5.4e-01& 10 \\  
8e+02&     1&     3&     3& 1&6.3e-01& 34 \\  
9e+02&     1&     3&     3& 1&6.4e-01& 31 \\  
1e+03&     1&     3&     3& 1&6.5e-01& 28 \\  
\belowspace
2e+03&     0&     0&     0& 1&1.0e+00& 6 \\ \hline
\end{tabular}
\end{sc}
%\end{small}
\end{center}
\vskip -0.1in
\end{table}

\end{document}